\documentclass{article}
\usepackage{amsmath}[1996/11/01]
\usepackage{amssymb,amsthm,amsxtra}
\usepackage{epic,eepic}
\usepackage{epsfig}

\def\[#1\]{\begin{equation}#1\end{equation}}
\makeatletter
\def\beq{%
   \relax\ifmmode
      \@badmath
   \else
      \ifvmode
         \nointerlineskip
         \makebox[.6\linewidth]%
      \fi
      $$
   \fi
}
\def\eeq{%
   \relax\ifmmode
      \ifinner
         \@badmath
      \else
         $$
      \fi
   \else
      \@badmath
   \fi
   \ignorespaces
}

\def\enddisplaymath{\eeq\global\@ignoretrue}
\makeatother

\newtheorem{thm}{Theorem}
\newtheorem{cor}[thm]{Corollary}
\newtheorem{lem}[thm]{Lemma}

\theoremstyle{remark}
\newtheorem*{rem}{Remark}
\newtheorem{rems}{Remark}[thm]

\theoremstyle{definition}
\newtheorem{defn}{Definition}

\numberwithin{equation}{section}
\numberwithin{thm}{section}

\DeclareMathOperator{\Tr}{Tr}
\DeclareMathOperator{\Exp}{E}
\DeclareMathOperator{\Aut}{Aut}
\renewcommand{\Re}{\operatorname{Re}}
\DeclareMathOperator{\sgn}{sgn}
\DeclareMathOperator{\Prob}{Pr}

\DeclareMathOperator{\pf}{pf}

\newcommand{\C}{\mathbb C}

\newcommand{\Z}{\mathbb Z}
\newcommand{\N}{\mathbb N}
\newcommand{\Mult}{\mathbb M}
\newcommand{\calM}{\mathcal M}
\newcommand{\defeq}{:=}

\newcommand\psymmU{%
\begin{picture}(1,1)(0,0)%
\allinethickness{0.5pt}%
\path(0,0)(0,1)(1,1)(1,0)(0,0)%
\end{picture}}
\newcommand\psymmUU{%
\begin{picture}(1,1)(0,0)%
\allinethickness{0.5pt}%
\path(0,0)(0,1)(1,1)(1,0)(0,0)%
\put(0.5,0.5){\makebox(0,0){$\cdot$}}%
\end{picture}}
\newcommand\psymmO{%
\begin{picture}(1,1)(0,0)%
\allinethickness{0.5pt}%
\path(0,0)(0,1)(1,1)(1,0)(0,0)%
\path(0,0)(1,1)%
\end{picture}}
\newcommand\psymmS{%
\begin{picture}(1,1)(0,0)%
\allinethickness{0.5pt}%
\path(0,0)(0,1)(1,1)(1,0)(0,0)%
\path(1,0)(0,1)%
\end{picture}}
\newcommand\psymmu{%
\begin{picture}(1,1)(0,0)%
\allinethickness{0.5pt}%
\path(0,0)(0,1)(1,1)(1,0)(0,0)%
\path(0,0)(1,1)%
\path(0,1)(1,0)%
\end{picture}}

\newbox\tsymmUbox
\newbox\tsymmUUbox
\newbox\tsymmObox
\newbox\tsymmSbox
\newbox\tsymmubox
\setbox\tsymmUbox =\hbox{\kern0.75pt\setlength{\unitlength}{6pt}\psymmU \kern0.75pt}
\setbox\tsymmUUbox=\hbox{\kern0.75pt\setlength{\unitlength}{6pt}\psymmUU\kern0.75pt}
\setbox\tsymmObox =\hbox{\kern0.75pt\setlength{\unitlength}{6pt}\psymmO \kern0.75pt}
\setbox\tsymmSbox =\hbox{\kern0.75pt\setlength{\unitlength}{6pt}\psymmS \kern0.75pt}
\setbox\tsymmubox =\hbox{\kern0.75pt\setlength{\unitlength}{6pt}\psymmu \kern0.75pt}
\def\tsymmU{{\copy\tsymmUbox}}
\def\tsymmUU{{\copy\tsymmUUbox}}
\def\tsymmO{{\copy\tsymmObox}}
\def\tsymmS{{\copy\tsymmSbox}}
\def\tsymmu{{\copy\tsymmubox}}
\newbox\symmUbox
\newbox\symmUUbox
\newbox\symmObox
\newbox\symmSbox
\newbox\symmubox
\setbox\symmUbox =\hbox{\kern0.75pt\setlength{\unitlength}{4.5pt}\psymmU \kern0.75pt}
\setbox\symmUUbox=\hbox{\kern0.75pt\setlength{\unitlength}{4.5pt}\psymmUU\kern0.75pt}
\setbox\symmObox =\hbox{\kern0.75pt\setlength{\unitlength}{4.5pt}\psymmO \kern0.75pt}
\setbox\symmSbox =\hbox{\kern0.75pt\setlength{\unitlength}{4.5pt}\psymmS \kern0.75pt}
\setbox\symmubox =\hbox{\kern0.75pt\setlength{\unitlength}{4.5pt}\psymmu \kern0.75pt}
\def\symmU{{\copy\symmUbox}}
\def\symmUU{{\copy\symmUUbox}}
\def\symmO{{\copy\symmObox}}
\def\symmS{{\copy\symmSbox}}
\def\symmu{{\copy\symmubox}}

\def\tsymmg{\circledast}
\def\symmg{\circledast}

\begin{document}

\title{{\bf Algebraic aspects of increasing subsequences}}
\author{{\bf Jinho Baik\footnote{Address after Sep. 1, 1999 : Princeton
University and Institite for Advanced Study, New Jersey}}\\Courant
Institute of Mathematical Sciences, New York\\baik@cims.nyu.edu\\
{\bf Eric M. Rains}\\AT\&T Research, New Jersey
\\rains@research.att.com}
\date{September 21, 2000}
\maketitle

\begin{abstract}
We present a number of results relating partial Cauchy-Littlewood sums,
integrals over the compact classical groups, and increasing subsequences of
permutations.  These include: integral formulae for the distribution of the
longest increasing subsequence of a random involution with constrained
number of fixed points; new formulae for partial Cauchy-Littlewood sums, as
well as new proofs of old formulae; relations of these expressions to
orthogonal polynomials on the unit circle; and explicit bases for invariant
spaces of the classical groups, together with appropriate generalizations
of the straightening algorithm.

\end{abstract}

\section*{Introduction}

Consider the following two identities:
\begin{align}
\sum_\lambda
s_\lambda(x_1,x_2,\ldots )
s_\lambda(y_1,y_2,\ldots )
&=
\prod_{i,j} (1-x_iy_j)^{-1}\label{eq:cauchys}\\
\lim_{l\to\infty}
\Exp_{U\in U(l)}
\prod_i \det(1-x_i U)^{-1}
\prod_i \det(1-y_i U^{-1})^{-1}
&=
\prod_{i,j} (1-x_iy_j)^{-1}.\label{eq:cauchyU}
\end{align}
The first of these is the well-known identity of Cauchy (\cite{Macdonald}).
The second is a formal analogue of the Szeg\"o limit theorem, equivalent to
a theorem of \cite{DiaconisShahshahani}.  Since the right-hand sides are
the same, we also have a third identity:
\[
\sum_\lambda
s_\lambda(x_1,x_2,\ldots )
s_\lambda(y_1,y_2,\ldots )
=
\lim_{l\to\infty}
\Exp_{U\in U(l)}
\prod_i \det(1-x_i U)^{-1}
\prod_i \det(1-y_i U^{-1})^{-1}.\label{eq:cauchysU}
\]
Our object of study in the present work is generalizations of these three
identities.

Our generalizations take two forms.  One is to remove the limit
$l\to\infty$.  As we shall see (section \ref{sec:symmfunc}), in order to
preserve \eqref{eq:cauchysU}, 
we must then restrict the sum over partitions
to involve only partitions with at most $l$ parts.  It is here that
increasing subsequences appear: in order to rescue equations
\eqref{eq:cauchys} and \eqref{eq:cauchyU}, we must replace the common
right-hand side with a generating function counting objects (reducing to
permutations in an appropriate limit) without long increasing subsequences.
In the case of \eqref{eq:cauchys}, the connection is via the
Robinson-Schensted-Knuth correspondence, with its well-known connections to
increasing subsequences.  It turns out that there is also a direct
connection for \eqref{eq:cauchyU}, in terms of the invariant theory of the
unitary group.  In particular, this gives a direct (and essentially
elementary) proof of the known connection between unitary group integrals
and increasing subsequences (\cite{Rains}), as well as of the new
connections given here.

The other way in which we generalize these identities is to replace the
unitary group $U(l)$ by one of four other groups, including the orthogonal
and symplectic groups.  In terms of permutations, this corresponds to
considering involutions (in two ways), signed permutations, and signed
involutions, in addition to the original case of permutations; each of
these conditions can be described as a particular type of symmetry
condition.  In each case, we obtain analogues of the finite $l$ versions of
equations \eqref{eq:cauchys}, \eqref{eq:cauchyU}, and \eqref{eq:cauchysU},
together with increasing subsequence interpretations.

\noindent {\bf Guide to main results}

One of the classical models used to analyze increasing subsequences is the
Poisson model: one generates a random subset of the unit square via a
Poisson process, then associates a permutation to this subset in a
canonical way (the order of the $y$ coordinates relative to the $x$
coordinates).  The five symmetry types correspond to the five subgroups of
$\Z_2\times \Z_2$ (acting on the unit square via diagonal reflections); one
insists that the random subset be preserved by the appropriate group.  It
turns out that each symmetry type is naturally associated with a certain
compact Lie group, determined as a fixed subgroup via an appropriate action
of $\Z_2\times \Z_2$ on the unitary group.  Our first main result (Theorem
\ref{thm:trmoments}) then states that the (exact) distribution of the
(length of the) longest increasing subsequence of a random permutation of a
given symmetry type is given by the moments of the trace of a random
element of the corresponding Lie group.

By the Schensted correspondence, each of the five cases of Theorem
\ref{thm:trmoments} can be viewed as expressing an integral over one of the
five groups as an appropriate sum over partitions.  Each such identity
specializes an appropriate Schur function identity (Theorem 5.2 and
Corollary 5.3).  For the three symmetry types with diagonal symmetry, this
can be further generalized (essentially allowing points on the diagonal of
symmetry); thus, for instance, if $f(\lambda)$ is the number of odd parts
and $\ell(\lambda)$ is the number of parts of $\lambda$, then (Theorem 5.6)
\[
\sum_{\ell(\lambda)\le l} \alpha^{f(\lambda)} s_\lambda(x)
=
E_{U\in O(l)} \det(1+\alpha U) \prod_j \det(1-x_j U)^{-1}.
\]

For a quite general class of specializations (corresponding to
``super-Schur'' or ``hook Schur'' functions), these Schur function sums
have combinatorial interpretations in terms of increasing subsequences of
multisets.  That is, for each symmetry type and each specialization
(subject to convergence conditions), we construct a random multiset and a
notion of increasing subsequence such that (Theorem 7.1) the normalized
Schur function sum gives the distribution of the longest increasing
subsequence.  These random multiset models generalize Johansson's 2D random
growth model \cite{kurtj:shape}.

Putting these together, we obtain a connection between integrals and
increasing subsequences of multisets.  In each case, the identity states
that the dimension of a certain space of invariants is given by counting a
certain collection of multisets without long increasing subsequences.  For
the three classical groups, we give direct proofs of this fact by producing
an explicit basis of the invariants indexed by the appropriate multisets.
Thus, for instance (Theorem 8.2), the centralizer algebra of the $n$th
tensor power representation of $U(l)$ has basis given by permutations of
length $n$ with no increasing subsequence of length greater than $l$.  More
generally (Theorem 8.4), the space of simultaneous multilinear invariants
of a collection of symmetric and antisymmetric, covariant and contravariant
tensors is explicitly indexed by multisets without long increasing
subsequences (generalizing the classical straightening algorithm).
Corresponding results hold for the orthogonal (Theorems 8.5 and 8.6) and
symplectic (Theorems 8.7 and 8.8) groups.

The remaining collection of results is of lesser interest combinatorially,
but is crucial to our asymptotic analysis in \cite{PartII}.  The key step
in the analysis is to express the integrals of interest in terms of
orthogonal polynomials on the unit circle.  This is done for the classical
groups in Theorem 2.3 (the remaining two groups reduce to the unitary
group); for the unitary group the connection is immediate, while for the
orthogonal group one must pass through orthogonal polynomials on $[-1,1]$.
We also give a number of results indicating how certain modifications to
the integrand affect the integral.  As a consequence, we find (Corollary
4.3) that for each of the five natural Poisson models, the longest
increasing subsequence distribution can be expressed in terms of the same
family of orthogonal polynomials.  In the sequel to this paper \cite{PartII},
we determine the asymptotics of such polynomials via the Riemann-Hilbert
method, and thus obtain the limiting longest increasing subsequence
distribution for each of the five Poisson models as well as for the
de-Poissonized versions (random symmetric permutations).  The results
are expressed in terms of the solution to the Painlev\'e II equation,
thus connecting to random matrix theory.  Further related asymptotic work
can be found in \cite{BR4}, \cite{Baik1}, \cite{Rains:mean_identity}.

\noindent {\bf Outline}

Section \ref{sec:symmiseq} introduces the five symmetry types and their
associated groups.  In particular, we express the (exact) distribution of
the longest increasing subsequence of a random permutation of a given
symmetry type as an integral over the corresponding group.
Section \ref{sec:dets1} expresses these integrals as determinants of Bessel
functions related to orthogonal polynomials; this relation is given in
section \ref{sec:orthopoly}.

Section \ref{sec:diags} describes the extension of the integrals to include
the cases when fixed points are allowed.  In order to prove the formulae of
section \ref{sec:diags}, section \ref{sec:symmfunc} uses
representation-theoretic arguments to deduce integral representations of
partial Cauchy-Littlewood sums, at which point the theory of symmetric
functions can be applied.  Section \ref{sec:pfaffians} briefly discusses
alternate proofs of those formulae based on intermediate pfaffian forms.

Section \ref{sec:moreiseq} uses a generalized Robinson-Schensted-Knuth
correspondence from \cite{BereleRemmel} to relate the partial
Cauchy-Littlewood sums to increasing subsequences of certain distributions
of random multisets.

Finally, section \ref{sec:invariants} shows that there is an intimate
connection between increasing subsequences and invariants of the classical
groups.  Indeed, all of the integrals for which we give increasing
subsequence interpretations also can be interpreted as dimensions of
certain spaces of invariants.  We give direct, elementary, proofs of these
identities, by constructing bases of invariants explicitly indexed by
multisets with restricted increasing subsequenc length.  In the process, we
obtain a generalization of the straightening algorithm of invariant theory,
as well as analogues for the orthogonal and symplectic groups.  We also
discuss extensions to quantum groups and supergroups.

\noindent {\bf Notation}
We refer the reader to \cite{Macdonald} for
notation and introduction to symmetric functions.  In other notation, if
$G$ is a compact group, we use
\[
\Exp_{U\in G} f(U)
\]
to denote the integral of $f(U)$ with respect to the (normalized) Haar
measure on $G$.  In other words, this is the expected value of $f$
evaluated at a uniform random element of $G$.  When $G$ is the orthogonal
group, we will occasionally need to consider the two components of
$G$.  Thus we will write
\[
\Exp_{U\in O^\pm(l)} f(U)
\]
to denote the integral of $f$ over the coset of $O(l)$ of determinant $\pm
1$.  In particular, we have:
\[
\Exp_{U\in O^\pm(l)} f(U)
=
{1\over 2} \Exp_{U\in O(l)} (1\pm \det(U)) f(U).
\]

\noindent {\bf Acknowledgments.}
We would like to acknowledge the following people for helpful discussions:
Kurt Johansson for telling us about the processes generalized in section
\ref{sec:moreiseq}, Richard Stanley for telling us about the references
for that section, Peter Shor for spotting flaws in earlier versions of the
algorithms of section \ref{sec:invariants}, and Christian Krattenthaler
for helpful comments on section \ref{sec:symmfunc}.  We would also like
to thank Jeff Lagarias, Andrew Odlyzko, and Neil Sloane for helpful
comments and enthusiasm.  The work of the first author was supported in
part by a Sloan Doctoral Foundation Fellowship.

\section{Symmetrized increasing subsequence problems}\label{sec:symmiseq}

One of the standard models used in analyzing the usual increasing
subsequence problem is defined as follows.  We say that a collection of $k$
points in the unit square is ``increasing'' if for any two points
$(x_1,y_1)$ and $(x_2,y_2)$, either $x_1<x_2$ and $y_1<y_2$ or $x_1>x_2$
and $y_1>y_2$.  One can ask, then, for the distribution of the size of the
largest increasing subset of $n$ points chosen uniformly and independently
from the unit square.  It is not too difficult to see that this is the same
as the distribution of the longest increasing subsequence of a random
permutation; indeed, we can associate a (uniformly distributed) permutation
to a given collection of points by using the relative order of the $y$
coordinates after sorting along the $x$ coordinates.  Also of interest
is the Poissonized analogue, in which new points are occasionally
added in such a way that the number of points at time $\lambda$ is
Poisson with parameter $\lambda$.

One way to generalize this model is to impose a symmetry condition on the
set of points.  The square has 8 symmetries; if we insist that the symmetry
preserve increasing collections, we obtain a group $H$ of 4 elements,
generated by the reflections through the main diagonals.  Thus there are 5
possible symmetry conditions we can impose (including the trivial
condition), which we will denote by the symbols $\tsymmU$, $\tsymmO$,
$\tsymmS$, $\tsymmUU$, and $\tsymmu$, with associated groups $H_\symmU$,
$H_\symmO$, $H_\symmS$, $H_\symmUU$, and $H_\symmu$.  (The symbol indicates
the point/line(s) of reflection) We will also use the symbol $\tsymmg$ to
denote an arbitrary choice of the five possibilities.

\begin{defn}
We define
\beq
p^\symmg_{nl}
\eeq
to be the probability that, if $n$ points are chosen uniformly at random in the
unit square, then the set $\Sigma$ consisting of the images of those points
under $H_\symmg$ contains no increasing subset with more than $l$
points.  We also define a function
\[
Q^\symmg_l(\lambda)=e^{-\lambda} \sum_{0\le n}
\frac{\lambda^n p^\symmg_{nl}}{n!}.
\]
\end{defn}

The function $Q^\symmg_l(\lambda)$ corresponds to the natural Poisson
model; $Q^\symmg_l(\lambda)$ is the probability that the largest increasing
subset at time $\lambda$ has size at most $l$, and $1-Q^\symmg_l(\lambda)$
is the distribution function of the time at which an increasing subset of
size $l+1$ first appears.  In the sequel, however, it will turn out to be
convenient to use a different time scale; we thus define
\[
P^\symmg_l(t)=Q^\symmg_l(a_\symmg t^2/2),
\]
where $a_\symmO=a_\symmS=1$, $a_\symmU=a_\symmu=2$, and $a_\symmUU=4$.
(Here $a_\symmg$ is the number of Young tableaux associated to an
element of $S^\symmg_n$ (q.v.))

If we map the set of points to a permutation, we obtain a permutation
uniformly distributed from an appropriate ensemble.  To be precise,
define the involution $\iota\in S_n$ by $x\mapsto n+1-x$.  Then
define an ensemble $S^\symmg_n$ for each symmetry type as follows:
{
\allowdisplaybreaks
\begin{align}
S^\symmU_n &= S_n\\
S^\symmO_n &= \{\pi\in S_{2n}:\pi=\pi^{-1},\ \pi(x)\ne x\}\\
S^\symmS_n &= \{\pi\in S_{2n}:\pi=\iota \pi^{-1}\iota ,\ \pi(x)\ne
\iota(x)\}\\
S^\symmUU_n &= \{\pi\in S_{2n}:\pi=\iota\pi\iota\}\\
S^\symmu_n &= \{\pi\in S_{4n}:\pi=\pi^{-1},\ \pi=\iota\pi^{-1}\iota,\ 
\pi(x)\ne x,\iota(x)\}.
\end{align}
}

\begin{lem}
If a set $\Sigma$ is chosen as above with symmetry $\tsymmg$,
then with probability 1, the associated permutation is well-defined
and is uniformly distributed from $S^\symmg_n$.
\end{lem}

This motivates the further definition
\[
f^\symmg_{nl} = |S^\symmg_n| p^\symmg_{nl}.
\]
That is, $f^\symmg_{nl}$ is the number of elements of
$S^\symmg_n$ with no increasing subsequence of length greater than
$l$.

It is straightforward to compute $|S^\symmg_n|$ for each case:
{
\allowdisplaybreaks
\begin{align}
|S^\symmU_n| &= n!\\
|S^\symmO_n| = |S^\symmS_n| &= {(2n)!\over 2^nn!}\\
|S^\symmUU_n| &= 2^nn!\\
|S^\symmu_n| &= {(2n)!\over n!}.
\end{align}
}
Note that
\[
P^\symmg_l(t) = e^{-a_\symmg t^2/2} \sum_{0\le n} f^\symmg_{nl}
{t^{2n}\over n!n!}
\]
for $\tsymmU$ and $\tsymmUU$, and
\[
P^\symmg_l(t) = e^{-a_\symmg t^2/2} \sum_{0\le n} f^\symmg_{nl}
{t^{2n}\over (2n)!}
\]
for $\tsymmO$, $\tsymmS$ and $\tsymmu$.

A major reason for considering the above problems is the following:

\begin{thm}\label{thm:trmoments}
Fix an integer $l>0$. Map $H$ into $\Aut(U(2l))$ by
\[
\diagup\mapsto (U\mapsto (U^t)^{-1})\quad\text{and}\quad
\diagdown\mapsto (U\mapsto -J(U^t)^{-1}J),
\]
where $J=\pmatrix 0&-I_l\\I_l&0\endpmatrix$.  Let $U^\symmg(2l)$
be the subgroup of $U(2l)$ fixed by the corresponding automorphisms.
Then
\[
f^\symmg_{n(2l)} = \Exp_{U\in U^\symmg(2l)} \left(|\Tr(U)|^{2n}\right).
\]
\end{thm}

Before giving the proof, it will be helpful to list the groups
$U^\symmg(2l)$:
{
\allowdisplaybreaks
\begin{align}
U^\symmU(2l) &= U(2l)\\
U^\symmO(2l) &= O(2l)\\
U^\symmS(2l) &= Sp(2l)\\
U^\symmUU(2l) &\cong U(l)\times U(l)\\
U^\symmu(2l) &= O(2l)\cap Sp(2l) \cong U(l).
\end{align}
}
The last instance is the image of the fundamental representation of
$U(l)$ as a $2l\times 2l$ real matrix, and thus corresponds to the direct
sum of the fundamental representation and its conjugate.

\begin{proof}
The cases $\tsymmU$, $\tsymmO$, and $\tsymmS$ are given in \cite{Rains}; more
precisely, $\tsymmU$ is given there as Theorem 1.1, while $\tsymmO$ and
$\tsymmS$ are given in Theorem 3.4.  (Note that if $\pi\in S^{\symmS}_n$,
then $\pi\iota$ is a fixed-point-free involution with decreasing
subsequences corresponding to the increasing subsequence of $\pi$.)
We also give new, elementary, proofs below (Theorems \ref{thm:invU},
\ref{thm:invO}, and \ref{thm:invS}).

It remains to consider $\tsymmUU$ and $\tsymmu$.
Via the Robinson-Schensted correspondence (see \cite{Kn}, section 5.1.4, for an
excellent introduction), we can associate a pair $(P,Q)$ of Young tableaux
of the same shape to $\pi\in S^{\symmUU}_n$, satisfying the relations
$P^S=P$, $Q^S=Q$, where $S$ is the duality operation (``evacuation'') of
Sch\"utzenberger (ibid.).  But there is a bijective correspondence between
self-dual tableaux and domino tableaux (see, e.g., \cite{vL96}), and
further from domino tableaux to pairs of ordinary tableaux with disjoint
content, and with shape determined only by the shape of the domino tableau
(\cite{StantonWhite}).  Thus we have associated four Young tableaux
$(P_1,P_2,Q_1,Q_2)$, where $P_1$ and $Q_1$ have the same shape, $P_2$,
$Q_2$ have the same shape, $P_1$ and $P_2$ have disjoint content, and $Q_1$
and $Q_2$ have disjoint content.  This corresponds to a choice of $0\le
m\le n$, independent choices of two subsets of size $m$ of $[1,2,\dots n]$,
and independent choices of $\pi_1$ and $\pi_2$ of length $m$ and $n-m$.
Furthermore, the longest increasing subsequence of $\pi$ has length at most
$2l$ precisely when the longest increasing subsequences of $\pi_1$ and
$\pi_2$ are of length at most $l$.  Putting this together, we find
\[
 f^\symmUU_{n(2l)} = \sum_{0\le m\le n} \binom{n}{m}^2
f^\symmU_{ml} f^\symmU_{(n-m)l}.
\]
On the other hand, the integral formula simplifies as follows:
\begin{align}
\Exp_{U\in U^\symmUU(2l)} \left(|\Tr(U)|^{2n}\right)
&=
\Exp_{U_1\in U(l),U_2\in U(l)} \left(|\Tr(U_1)+\Tr(U_2)|^{2n}\right)\\
&=
\Exp_{U_1\in U(l),U_2\in U(l)}
\left(
|(\Tr(U_1)+\Tr(U_2))^n|^2\right)\\
&=
\Exp_{U_1\in U(l),U_2\in U(l)}
\Big(
\Big|\sum_{0\le m\le n} \binom{n}{m} \Tr(U_1)^m\Tr(U_2)^{n-m}\Big|^2\Big)\\
&=
\sum_{0\le m\le n}
\binom{n}{m}^2
f^\symmU_{ml}f^\symmU_{(n-m)l}\\
&=
f^\symmUU_{n(2l)}.
\end{align}

Similarly, an element $\pi\in S^\symmu_n$ corresponds to a pair of
Young tableaux of the same shape with disjoint content, and thus
\[
f^\symmu_{n(2l)}=\binom{2n}{n} f^\symmU_{nl}.
\]
On the other hand,
\begin{align}
\Exp_{U\in U^\symmu(2l)} \left(|\Tr(U)|^{2n}\right)
&=
\Exp_{U\in U(l)} \left((\Tr(U)+\overline{\Tr(U)})^{2n}\right)\\
&=
\sum_{0\le m\le 2n}
\binom{2n}{m}
\Exp_{U\in U(l)} \left(\Tr(U)^m\overline{\Tr(U)}^{2n-m}\right)\\
&=
\binom{2n}{n}
\Exp_{U\in U(l)} \left(|\Tr(U)|^{2n}\right).
\end{align}
\end{proof}

It ought to be possible to give a more uniform proof of this result; the
results of Section \ref{sec:invariants} (q.v.) may be relevant to this goal.

There is an analogue of Theorem \ref{thm:trmoments} in which $2l$ is
replaced by $2l+1$:

\begin{thm}
For any $n$, $l\ge 0$,
\begin{align}
f^\symmU_{n(2l+1)} &= \Exp_{U\in U(2l+1)}           |\Tr(U)|^{2n}\\
f^\symmO_{n(2l+1)} &= \Exp_{U\in O(2l+1)}           |\Tr(U)|^{2n}\\
f^\symmUU_{n(2l+1)}&= \Exp_{U\in U(l)\oplus U(l+1)} |\Tr(U)|^{2n},
\end{align}
while
\begin{align}
f^\symmS_{n(2l+1)}&=f^\symmS_{n(2l)}\\
f^\symmu_{n(2l+1)}&=f^\symmu_{n(2l)}.
\end{align}
\end{thm}

Also, we have the following corollary for $\tsymmUU$ and $\tsymmu$:

\begin{cor}
For any $n$, $l\ge 0$,
{
\allowdisplaybreaks
\begin{align}
P^\symmUU_{2l}(t)&=P^\symmU_l(t)^2\\
P^\symmUU_{2l+1}(t)&=P^\symmU_l(t)P^\symmU_{l+1}(t)\\
P^\symmu_{2l}(t)&=P^\symmU_l(t).
\end{align}
}
\end{cor}

And for $\tsymmO$, $\tsymmS$, and $\tsymmu$, we have

\begin{cor}\label{cor:oldcor6}
For any $n$, $l\ge 0$,
{
\allowdisplaybreaks
\begin{align}
P^\symmO_l(t)&=e^{-t^2/2} \Exp_{U\in O(l)} \exp(t\Tr(U))\\
P^\symmS_{2l}(t)&=e^{-t^2/2} \Exp_{U\in Sp(2l)} \exp(t\Tr(U))\\
P^\symmu_{2l}(t)&=e^{-t^2} \Exp_{U\in U^\symmu(2l)} \exp(t\Tr(U)).
\label{eq:opww}
\end{align}
}
\end{cor}

\begin{proof}
For $\tsymmO$, we have
\[
e^{t^2/2} P^\symmO_{l}(t)
=
\sum_{0\le n}
{t^{2n}\over (2n)!}
\Exp_{U\in O(l)} \left(\Tr(U)^{2n}\right)
\]
But $\Exp_{U\in O(l)} \Tr(U)^n=0$ for $n$ odd, so this is
\[
\Exp_{U\in O(l)} \exp(t\Tr(U)),
\]
as required.  The calculations for $\tsymmS$ and $\tsymmu$ are analogous.
\end{proof}

\begin{rem}
In particular, we see that the formula \eqref{eq:opww} which was derived
in \cite{OPWW} as an expression for $\tsymmU$ is really most naturally
interpreted in $\tsymmu$ terms.
\end{rem}

As an aside, we observe that if we remove the condition that the symmetries
under consideration preserve increasing sets, but insist that the
corresponding sets should still give permutations, there is one further
type of symmetry allowed, namely rotation by 90 degrees.  In terms of
permutations, this is the set
\[
S^\circ_n=\{\pi\in S_{4n}:\pi^2=\iota\}\qquad |S^\circ_n|=(2n)!/n!
\]
Such permutations correspond to pairs of tableaux $(P,Q^t)$ with $n$
elements such that $P$ and $Q$ have the same shape and disjoint content.
It follows that the length $l$ of the longest increasing subsequence from
this set has the same distribution as $\max(2 l^+(\pi),2 l^-(\pi)-1)$,
where $\pi$ is randomly chosen from $S_n$, $l^+(\pi)$ is the increasing
subsequence length of $\pi$, and $l^-(\pi)$ is the decreasing subsequence
length of $\pi$.  In particular, the bound of Erdos and Szekeres implies
that $f^\circ_{nl}=0$ for $n>l^2$, and thus no integral formula \`a la
Theorem \ref{thm:trmoments} can exist for this case.  There is a
determinant formula, however, which can be obtained from the following
symmetric function identity:
\[
\sum_{\substack{\ell(\lambda)\le l^+\\\ell(\lambda')\le l^-}} s_\lambda(x) s_{\lambda'}(y)
=
\det
\pmatrix (h_{j - i}(x))_{0\le i<l^+, 0\le j<l^+ + l^-}\\
      ((-1)^{l^+ + i-j} h_{l^+ +i-j}(y))_{0\le i<l^-, 0\le j<l^+ + l^-} )
\endpmatrix.
\]
One can either derive this via the approach in \cite{Goulden}, or simply
use the Jacobi-Trudi identity together with matrix manipulation as in
\cite{Gessel}.  This then gives the following formulae:
\begin{align}
\sum_{0\le n} f^\circ_{n(2l)} {t^{2n}\over (2n)!}
&=
\det\pmatrix
({t^{j - i}\over (j-1)!})_{0\le i<l, 0\le j<2l}\\
({(-t)^{l + i-j}\over (l +i-j)!})_{0\le i<l, 0\le j<2l}
\endpmatrix\\
\sum_{0\le n} f^\circ_{n(2l+1)} {t^{2n}\over (2n)!}
&=
\det\pmatrix({t^{j - i}\over (j-1)!})_{0\le i<l, 0\le j<2l+1}\\
      ({(-t)^{l + i-j}\over (l +i-j)!})_{0\le i<l+1, 0\le j<2l+1}\endpmatrix.
\end{align}
It is not clear how to obtain asymptotic information from the formulae,
however, so we will not discuss this case further.  (The techniques of
\cite{BOO} may be applicable, however.)

\section{Determinantal formulae}\label{sec:dets1}


It turns out that each of the formulae of Corollary \ref{cor:oldcor6} can
be expressed in terms of a Toeplitz or Hankel determinant, related to
orthogonal polynomials on the unit circle.  Indeed, this correspondence is
more general.

For the unitary group:

\begin{thm}\label{thm2.1}
Let $f(z)$, $g(z)$ be any functions on the unit circle.  Then
\[
\Exp_{U\in U(l)} \det(f(U)) \det(g(U^\dagger))
=
\det\left(
{1\over 2\pi}
\int_{[0,2\pi]} f(e^{i\theta}) g(e^{-i\theta}) e^{i(j-k)\theta} d\theta
\right)_{0\le j,k<l}.
\]
\end{thm}

\begin{proof}
Using the Weyl integration formula for the unitary group, we have
\[
\Exp_{U\in U(l)} \det(f(U)) \det(g(U^\dagger))
=
{1\over l!(2\pi)^l}
\int_{[0,2\pi]^l} \prod_{0\le j<k<l} \bigl| e^{ij\theta}-e^{ik\theta}\bigr|^2
\prod_{0\le j<l}f(e^{i\theta_j}) g(e^{-i\theta_j}) d\theta_j.
\]

The result follows from the standard theory of Toeplitz determinants, or by
the classic formula for the integral of a product of two (generalized)
Vandermonde determinants (see, for instance, \cite{deBruijn}):
\[
{1\over l!}
\int_{S^l}
\!
\det(\phi_j(x_k))_{0\le j,k<l}
\det(\psi_j(x_k))_{0\le j,k<l}
\prod_j \mu(dx_j)
=
\det\left(\int_S \phi_j(x) \psi_k(x) \mu(dx)\right)_{0\le j,k<l},
\]
for any measure $\mu$ on any set $S$.
\end{proof}

Similarly, for the orthogonal and symplectic groups, we have:

\begin{thm}\label{thm:intstodets}
Let $g(z)$ be any function on the unit circle such that the integrals
\[
\iota_j = {1\over 2\pi} \int_{[0,2\pi]} g(e^{i\theta})g(e^{-i\theta})
e^{ij\theta} d\theta
\]
are well-defined.  Then
{
\allowdisplaybreaks
\begin{align}
\Exp_{U\in O^+(2l)} \det(g(U))
&=
{1\over 2}
\det(\iota_{j-k}+\iota_{j+k})_{0\le j,k<l}\\
\Exp_{U\in O^-(2l)} \det(g(U))
&=
g(1)g(-1)
\det(\iota_{j-k}-\iota_{j+k+2})_{0\le j,k<l-1}\\
\Exp_{U\in O^+(2l+1)} \det(g(U))
&=
g(1)
\det(\iota_{j-k}-\iota_{j+k+1})_{0\le j,k<l}\\
\Exp_{U\in O^-(2l+1)} \det(g(U))
&=
g(-1)
\det(\iota_{j-k}+\iota_{j+k+1})_{0\le j,k<l}\\
\Exp_{U\in Sp(2l)} \det(g(U))
&=
\det(\iota_{j-k}-\iota_{j+k+2})_{0\le j,k<l},
\end{align}
}
except that $\Exp_{U\in O^+(0)} \det(g(U))=1$.
\end{thm}

\begin{proof}
As observed in Proposition 3.1 of \cite{Jo1}, integrals over the orthogonal and
symplectic groups can be expressed as Hankel determinants; thus,
for instance,
\[
\Exp_{U\in O^+(2l)} \det(g(U))
\propto
\det\left({1\over \pi} \int_{[0,2\pi]} g(e^{i\theta})g(e^{-i\theta})
p_j(\cos(\theta))p_k(\cos(\theta)) d\theta\right)_{0\le j,k<l}
\]
for any polynomials $p_j$ with $\deg(p_j)=j$.  In particular, this must be
true when $p_j(\cos(\theta))=\cos(j\theta)$ (Chebyshev polynomials). In
that case, 
noting that $\iota_k=\iota_{-k}$, the $jk$ coefficient of the determinant is
\[
{1\over \pi} \int_{[0,2\pi]} g(e^{i\theta})g(e^{-i\theta})
\cos(j\theta)\cos(k\theta) d\theta
=
\iota_{j-k}+\iota_{j+k}.
\]
The constant of proportionality can be determined by comparing the two
sides when $g=1$, and thus $\iota_j = \delta_{j0}$.

For the other cases, we take $p_j(\cos(\theta))$ to be
\[
{\sin((j+1)\theta)\over \sin(\theta)},\ 
{\sin((j+1/2)\theta)\over \sin(\theta/2)},\ 
{\cos((j+1/2)\theta)\over \cos(\theta/2)},\ 
\]
as appropriate.
\end{proof}

For our purposes, we will need the following related result:

\begin{thm}\label{thm:intstopolys}
Let $g(z)$ be as above, and consider the (symmetric) inner product on
polynomials
\[
\langle p(z),q(z)\rangle = {1\over 2\pi} \int_{[0,2\pi]}
p(e^{i\theta}) q(e^{-i\theta}) g(e^{i\theta}) g(e^{-i\theta}) d\theta.
\]
Let $\pi_j(z)$ be the monic orthogonal polynomials on the unit circle
relative to that inner product, and define $N_j=\langle
\pi_j(z),\pi_j(z)\rangle$.  If the polynomials $\pi_j(z)$ are well-defined,
then
{
\allowdisplaybreaks
\begin{align}
\Exp_{U\in U(l)} |\det(g(U))|^2
&=
\prod_{0\le j<l} N_j\\
\Exp_{U\in O^+(0)} \det(g(U))&= 1\\
\Exp_{U\in O^+(2l)} \det(g(U))
&=
N_0 \prod_{0\le j<l-1} N_{2j+2}(1+\pi_{2j+2}(0))^{-1}\\
\Exp_{U\in O^-(2l)} \det(g(U))
&=
g(1)g(-1)\prod_{0\le j<l-1} N_{2j+2}(1-\pi_{2j+2}(0))^{-1}\\
\Exp_{U\in O^+(2l+1)} \det(g(U))
&=
g(1)\prod_{0\le j<l} N_{2j+1}(1-\pi_{2j+1}(0))^{-1}\\
\Exp_{U\in O^-(2l+1)} \det(g(U))
&=
g(-1)\prod_{0\le j<l} N_{2j+1}(1+\pi_{2j+1}(0))^{-1}\\
\Exp_{U\in Sp(2l)} \det(g(U))
&=
\prod_{0\le j<l} N_{2j+2}(1-\pi_{2j+2}(0))^{-1}.
\end{align}
}
\end{thm}

The proof will be given in Section~\ref{sec:orthopoly} below.  By combining
these formulae, we obtain:

\begin{cor}
Let $g(z)$ be as above.  Then for any $l>0$,
\[
g(1)g(-1) \Exp_{U\in U(l)} |\det(g(U))|^2
=
(\Exp_{U\in O^+(l+1)} \det(g(U)))(\Exp_{U\in O^-(l+1)} \det(g(U))).
\]
\end{cor}

\begin{proof}
Use $N_l=(1-\pi_l(0)^2)N_{l-1}$ which follows from \eqref{8.11} below.
\end{proof}

In our case, the function $g(z)=e^{tz}$, and thus everything
is related to the orthogonal polynomials on the circle for the weight
function $\exp(t(z+1/z))=\exp(2t\cos\theta)$.  Let these polynomials
be $\pi_l(z;t)$ with norms $N_l(t)$, and define five functions
$D_l(t)$, $D^{\pm\pm}(t)$ by
\begin{align}
D_l(t)
&=\det(I_{j-k}(2t))_{0\le j,k<l}\\
&=\prod_{0\le j<l} N_j(t)\\
\displaybreak[0]
D^{--}_0(t)&=1\\
D^{--}_l(t)
&=
{1\over 2}\det(I_{j-k}(2t)+I_{j+k}(2t))_{0\le j,k<l}\\
&=
N_0(t) \prod_{0\le j<l-1} N_{2j+2}(t)(1+\pi_{2j+2}(0;t))^{-1}\\
\displaybreak[0]
D^{++}_l(t)&=
\det(I_{j-k}(2t)-I_{j+k+2}(2t))_{0\le j,k<l}\\
&=
\prod_{0\le j<l} N_{2j+2}(t)(1-\pi_{2j+2}(0;t))^{-1}\\
\displaybreak[0]
D^{+-}_l(t)&=
\det(I_{j-k}(2t)-I_{j+k+1}(2t))_{0\le j,k<l}\\
&=
\prod_{0\le j<l} N_{2j+1}(t) (1-\pi_{2j+1}(0;t))^{-1}\\
\displaybreak[0]
D^{-+}_l(t)&=
\det(I_{j-k}(2t)+I_{j+k+1}(2t))_{0\le j,k<l}\\
&=
\prod_{0\le j<l} N_{2j+1}(t) (1+\pi_{2j+1}(0;t))^{-1},
\end{align}
where
\[
I_j(2t)
={1\over 2\pi}\int_{[0,2\pi]} e^{2t\cos\theta} e^{ij\theta}
=
\sum_m {t^m\over m!}{t^{j+m}\over (j+m)!}
\]
is a modified Bessel function of the first kind.

Then
{
\allowdisplaybreaks
\begin{align}
\Exp_{U\in U(l)} |\exp(t\Tr(U))|^2 &= D_l(t)\\
\Exp_{U\in O^+(2l)} \exp(t\Tr(U)) &= D^{--}_l(t)\\
\Exp_{U\in O^-(2l)} \exp(t\Tr(U)) &= D^{++}_{l-1}(t)\\
\Exp_{U\in O^+(2l+1)} \exp(t\Tr(U)) &= e^t D^{+-}_l(t)\\
\Exp_{U\in O^-(2l+1)} \exp(t\Tr(U)) &= e^{-t} D^{-+}_l(t).
\end{align}
}

In summary:

\begin{thm}\label{thm2.5}
For $l\ge 0$, we have the following formulae:
{
\allowdisplaybreaks
\begin{align}
P^\symmU_l(t) &= e^{-t^2} D_l(t)\\
P^\symmO_{2l}(t) &= e^{-t^2/2}
\left[
D^{--}_l(t)+
D^{++}_{l-1}(t)\right]/2\\
P^\symmO_{2l+1}(t) &= e^{-t^2/2}
\left[
e^tD^{+-}_l(t)+
e^{-t}D^{-+}_l(t)\right]/2\\
P^\symmS_{2l}(t) &=
e^{-t^2/2} D^{++}_l(t),\\
P^\symmUU_{2l}(t) &= e^{-2t^2} D_l(t)^2\\
P^\symmUU_{2l+1}(t) &= e^{-2t^2} D_l(t)D_{l+1}(t)\\
P^\symmu_{2l}(t) &= e^{-t^2} D_l(t)
\end{align}
}
except that
$
P^\symmO_0(t) = e^{-t^2/2} D^{--}_0(t)
= e^{-t^2/2}
$.
\end{thm}

In \cite{PartII}, we will also need the following limits, which follow
immediately from Szeg\"o's limit theorem and the analogue \cite{Jo1} for
orthogonal polynomials on a finite interval:

\begin{thm}
For any real $t\ge 0$,
{
\allowdisplaybreaks
\begin{align}
\lim_{l\to\infty} D_l(t) &= e^{t^2}\\
\lim_{l\to\infty} D^{--}_l(t) &= e^{t^2/2}\\
\lim_{l\to\infty} D^{++}_l(t) &= e^{t^2/2}\\
\lim_{l\to\infty} D^{+-}_l(t) &= e^{-t+t^2/2}\\
\lim_{l\to\infty} D^{-+}_l(t) &= e^{t+t^2/2}.
\end{align}
}
\end{thm}

\begin{rem}
These limits are also valid as limits of formal power series in $t$;
see \cite{DiaconisShahshahani} and Theorem \ref{thm:formalSzego} below.
\end{rem}

\begin{cor}\label{cor1}
For any real $t\ge 0$,
\begin{eqnarray}
e^{-t^2}D_l(t) &=&\prod_{j\geq l} N_j(t)^{-1} \\
e^{-t^2/2}D^{--}_l(t) &=& 
\prod_{j\geq l} N_{2j+2}(t)^{-1}(1+\pi_{2j+2}(0;t))\\
e^{-t^2/2}D^{++}_l(t) &=& 
\prod_{j\geq l} N_{2j+2}(t)^{-1}(1-\pi_{2j+2}(0;t)) \\
e^{-t^2/2+t}D^{+-}_l(t) &=& 
\prod_{j\geq l} N_{2j+1}(t)^{-1}(1-\pi_{2j+1}(0;t))\\
e^{-t^2/2-t}D^{-+}_l(t) &=& 
\prod_{j\geq l} N_{2j+1}(t)^{-1}(1+\pi_{2j+1}(0;t)).
\end{eqnarray}
\end{cor}

\section{Orthogonal polynomial identities}\label{sec:orthopoly}

Let $w(x)$ on $[-1,1]$ and $f(z)$ on the unit circle
be related by
\[
f(e^{i\theta})=w(\cos\theta).
\]
Associated to $f$ are five sets of monic orthogonal polynomials:
$\pi_l(z)$ with respect to the (symmetric) inner product
\[
\langle p(z),q(z)\rangle = {1\over 2\pi} \int_{[0,2\pi]} p(e^{i\theta})
q(e^{-i\theta}) f(e^{i\theta})\, d\theta.
\]
on the unit circle, and the four sets $p^{\pm\pm}_l(z)$, with respect to
the inner products
\[
\langle p(x),q(x)\rangle=
{1\over \pi} \int_{[-1,1]} p(x) q(x) w(x) (1-x)^{\pm 1/2}(1+x)^{\pm 1/2}\, dx.
\]
The notation $\langle p,q\rangle$ will refer to the inner product with respect
to $f(z)$ or $w(x)(1-x^2)^{-1/2}$, whichever is appropriate.
Thus the defining identities for these polynomials are:
{
\allowdisplaybreaks
\begin{align}
\langle \pi_n,\pi_m\rangle       &= \delta_{nm}N_n\\
\langle p^{--}_n,p^{--}_m\rangle &= \delta_{nm}N^{--}_n\\
\langle (1+x)p^{-+}_n,p^{-+}_m\rangle &= \delta_{nm}N^{-+}_n\\
\langle (1-x)p^{+-}_n,p^{+-}_m\rangle &= \delta_{nm}N^{+-}_n\\
\langle (1-x^2)p^{++}_n,p^{++}_m\rangle &= \delta_{nm}N^{++}_n.
\end{align}
}
We also use the notation
\[
\pi^*_l(z)=z^l \pi_l(1/z).
\]

Then the following identities hold (see for example Theorem 11.5 of
\cite{Szego}):
For $\pi$: 
{
\allowdisplaybreaks
\begin{align}
\pi_{l+1}(z)&=z\pi_l(z)+\pi_{l+1}(0)\pi^*_l(z)\\
\label{8.11} N_l &= N_0 \prod_{1\le j\le l} (1-\pi_j(0)^2)\\
\pi_l(1) &= \prod_{1\le j\le l} (1+\pi_j(0))\\
(-1)^l \pi_l(-1) &= \prod_{1\le j\le l} (1+(-1)^j\pi_j(0)).
\end{align}
}
For $p^{--}$:
{
\begin{align}
p^{--}_0((z+1/z)/2) &= 1\\
\displaybreak[0]
(2z)^l p^{--}_l((z+1/z)/2) &= \pi^*_{2l-1}(z)+z\pi_{2l-1}(z)\\
\displaybreak[0]
p^{--}_l(1)&=2^{-l}\prod_{0\le j<2l} (1+\pi_j(0))\\
\displaybreak[0]
p^{--}_l(-1)&=(-2)^{-l}\prod_{0\le j<2l} (1+(-1)^j \pi_j(0))\\
\displaybreak[0]
N^{--}_l
&=
4^{1/2-l}(1+\pi_{2l}(0))^{-1} \prod_{0\le j<2l} (1-\pi_j(0)^2)\\
&=
4^{1/2-l}N_{2l}(1+\pi_{2l}(0))^{-1}
\end{align}
}
For $p^{+-}$:
{
\begin{align}
(1-z) (2z)^l p^{+-}_l((z+1/z)/2)
&=
\pi^*_{2l}(z)-z\pi_{2l}(z)\\
\displaybreak[0]
p^{+-}_l(-1)&=(-2)^{-l}\prod_{1\le j\le 2l} (1+(-1)^j\pi_j(0))\\
\displaybreak[0]
N^{+-}_l
&=
4^{-l}(1+\pi_{2l+1}(0))
\prod_{1\le j\le 2l} (1-\pi_j(0)^2)\\
&=
4^{-l} N_{2l} (1+\pi_{2l+1}(0))\\
&=
4^{-l} N_{2l+1} (1-\pi_{2l+1}(0))^{-1}\\
\end{align}
}
For $p^{-+}$:
{
\begin{align}
(1+z) (2z)^l p^{-+}_l((z+1/z)/2)
&=
\pi^*_{2l}(z)+z\pi_{2l}(z)\\
\displaybreak[0]
p^{-+}_l(1)&=2^{-l}\prod_{1\le i\le 2l}(1+\pi_j(0))\\
\displaybreak[0]
N^{-+}_l
&=
4^{-l}(1-\pi_{2l+1}(0))
\prod_{1\le i\le 2l} (1-\pi_{i}(0)^2)\\
&=
4^{-l} N_{2l} (1-\pi_{2l+1}(0))\\
&=
4^{-l} N_{2l+1} (1+\pi_{2l+1}(0))^{-1}
\end{align}
}
And finally, for $p^{++}$:
{
\begin{align}
(1-z^2)(2z)^l p^{++}_l((z+1/z)/2)
&=
\pi^*_{2l+1}(z)-z\pi_{2l+1}(z)\\
\displaybreak[0]
N^{++}_l
&=
4^{-l-1/2}(1+\pi_{2l+2}(0)) \prod_{1\le j\le 2l+1} (1-\pi_j(0)^2)\\
&=
4^{-l-1/2} N_{2l+1}(1+\pi_{2l+2}(0))\\
&=
4^{-l-1/2} N_{2l+2}(1-\pi_{2l+2}(0))^{-1}
\end{align}
}

All of the proofs are straightforward.

We can now prove Theorem \ref{thm:intstopolys}.  Recall that we have chosen
$f$ to be of the form $f(z)=g(z)g(1/z)$ such that the above inner product
is well-defined and nondegenerate.  Now consider the integral for
$O^+(2l)$.  By the proof of Theorem \ref{thm:intstodets} and the theory of
Hankel determinants, we have
\[
\Exp_{U\in O^+(2l)} \det(g(U)) \propto \prod_{0\le j<l} N^{++}_j;
\]
but then this is in turn proportional to
\[
N_0 \prod_{0\le j<l-1} N_{2j+2}(1+\pi_{2j+2}(0))^{-1}.
\]
The constant can be determined by taking $g=1$.
The other cases are analogous.

\begin{thm}\label{thm:alphaformulae}
Let $f$ and $g$ be as above.  Then for any $\alpha$,
{
\allowdisplaybreaks
\begin{align}
\Exp_{U\in O^{\pm}(2l)} \det((1-\alpha U)g(U))
&=
(\pi^*_{2l-1}(\alpha)\pm\alpha\pi_{2l-1}(\alpha))
\Exp_{U\in O^{\pm}(2l)} \det(g(U))\\
\Exp_{U\in O^{\pm}(2l+1)} \det((1-\alpha U)g(U))
&=
(\pi^*_{2l}(\alpha)\mp\alpha\pi_{2l}(\alpha))
\Exp_{U\in O^{\pm}(2l+1)}\det(g(U))
\end{align}
}
As a special case,
{
\allowdisplaybreaks
\begin{align}
\Exp_{U\in O^+(l)} 
\det((1+U)g(U))
&=
2 g(-1)^{-1}
\Exp_{U\in O^-(l+1)} \det(g(U))
\\
\Exp_{U\in O^-(l)} 
\det((1+U)g(U))
&=
0\\
\Exp_{U\in O^+(2l)}
\det((1-U)g(U))
&=
2 g(1)^{-1}
\Exp_{U\in O^+(2l+1)} \det(g(U))
\\
\Exp_{U\in O^-(2l)}
\det((1-U)g(U))
&=
0
\\
\Exp_{U\in O^+(2l+1)} 
\det((1-U)g(U))
&=
0\\
\Exp_{U\in O^-(2l+1)} 
\det((1-U)g(U))
&=
2 g(1)^{-1}
\Exp_{U\in O^-(2l+2)} \det(g(U))
\end{align}
}
\end{thm}

\begin{proof}
Follows by standard results about the behavior of orthogonal polynomials
when the weight function is multiplied by a polynomial.
\end{proof}

Similarly,

\begin{thm}
With notation as above, and for $\alpha\beta\ne 1$,
\[
\Exp_{U\in U(l)} \det((1-\alpha U)(1-\beta U^\dagger)g(U)g(U^\dagger))
=
{\pi^*_l(\alpha)\pi^*_l(\beta)-\alpha\beta\pi_l(\alpha)\pi_l(\beta)
\over 1-\alpha\beta} \Exp_{U\in U(l)} \det(g(U)g(U^\dagger))
\]
In particular,
{
\begin{align}
2 g(1)g(-1)(1-\alpha\beta)\Exp_{U\in U(l)} &\det((1-\alpha U)(1-\beta U^\dagger)g(U)g(U^\dagger))\notag\\
=&
\Exp_{U\in O^+(l+1)} \det((1-\alpha U)g(U))
\Exp_{U\in O^-(l+1)} \det((1-\beta U)g(U))\\
+&
\Exp_{U\in O^+(l+1)} \det((1-\beta U)g(U))
\Exp_{U\in O^-(l+1)} \det((1-\alpha U)g(U))\notag
\end{align}
}
\end{thm}

For our purposes, the following form is preferable:

\begin{cor}
With notation as above, and for $\alpha\beta\ne 1$, $|\beta|<1$,
{
\begin{align}
{2 g(1)g(-1)(1-\alpha\beta)\over (1-\beta^2)}
\Exp_{U\in U(l)} &\det((1-\alpha U)(1-\beta U)^{-1}g(U)g(U^\dagger))\notag\\
=&
\Exp_{U\in O^+(l+1)} \det((1-\alpha U)(1-\beta U)^{-1}g(U))
\Exp_{U\in O^-(l+1)} \det(g(U))\\
+&
\Exp_{U\in O^+(l+1)} \det(g(U))
\Exp_{U\in O^-(l+1)} \det((1-\alpha U)(1-\beta U)^{-1}g(U))\notag
\end{align}
}
\end{cor}

\begin{rem}
This, together with the following theorem, enables us to compute the limits
$\alpha\to \pm 1$ or $\beta\to \pm 1$.
\end{rem}

\begin{thm}
With notation as above, and for $|\beta|<1$,
{
\allowdisplaybreaks
\begin{align}
\Exp_{U\in O^+(2l+2)} \det((1-\beta U)^{-1} g(U))/\Exp_{U\in O^+(2l+2)} \det(g(U))
&=
\frac{
{1\over 2}
\beta^{-l} \langle 2^l p^{--}_l(x),(1-\beta^2)/(1-2\beta x+\beta^2)\rangle
}{
(1-\beta^2) N_{2l} (1+\pi_{2l}(0))^{-1}
}
\\
\Exp_{U\in O^-(2l+2)} \det((1-\beta U)^{-1} g(U))/\Exp_{U\in O^-(2l+2)} \det(g(U))
&=
\frac{
\beta^{1-l} \langle 2^l (1-x^2) p^{++}_{l-1}(x),1/(1-2\beta x+\beta^2)
\rangle
}{
(1-\beta^2) N_{2l} (1-\pi_{2l}(0))^{-1}
}
\\
\Exp_{U\in O^+(2l+3)} \det((1-\beta U)^{-1} g(U))/\Exp_{U\in O^+(2l+3)} \det(g(U))
&=
\frac{
\beta^{-l} \langle 2^l (1-x) p^{+-}_l(x),(1+\beta)/(1-2\beta x+\beta^2)\rangle
}{
(1-\beta^2) N_{2l} (1+\pi_{2l+1}(0))
}
\\
\Exp_{U\in O^-(2l+3)} \det((1-\beta U)^{-1} g(U))/\Exp_{U\in O^-(2l+3)} \det(g(U))
&=
\frac{
\beta^{-l} \langle 2^l (1+x) p^{-+}_l(x),(1-\beta)/(1-2\beta x+\beta^2)\rangle
}{
(1-\beta^2) N_{2l} (1-\pi_{2l+1}(0))
}
\\
\Exp_{U\in Sp(2l)} \det((1-\beta U)^{-1} g(U))/\Exp_{U\in Sp(2l)} \det(g(U))
&=
\frac{
\beta^{1-l} \langle 2^l (1-x^2) p^{++}_{l-1}(x),1/(1-2\beta x+\beta^2)\rangle
}{
N_{2l} (1-\pi_{2l}(0))^{-1}
}
.
\end{align}
}
In the limits as $\beta\to \pm 1$,
{
\allowdisplaybreaks
\begin{align}
\lim_{\substack{\beta\to 1\\ |\beta|<1}}
(1-\beta^2) \Exp_{U\in O^{\pm}(l+1)} \det((1-\beta U)^{-1} g(U))
&=
g(1) \Exp_{U\in O^{\pm}(l)} \det(g(U))\\
\lim_{\substack{\beta\to -1\\ |\beta|<1}}
(1-\beta^2) \Exp_{U\in O^{\pm}(l+1)} \det((1-\beta U)^{-1} g(U))
&=
g(-1) \Exp_{U\in O^{\mp}(l)} \det(g(U))
\end{align}
}
For the symplectic group,
{
\allowdisplaybreaks
\begin{align}
\lim_{\substack{\beta\to 1\\ |\beta|<1}}
\Exp_{U\in Sp(2l)} \det((1-\beta U)^{-1} g(U))
&=
g(-1)^{-1} \Exp_{U\in O^-(2l+1)} \det(g(U))\\
\lim_{\substack{\beta\to -1\\ |\beta|<1}}
\Exp_{U\in Sp(2l)} \det((1-\beta U)^{-1} g(U))
&=
g(1)^{-1} \Exp_{U\in O^+(2l+1)} \det(g(U))
\end{align}
}
\end{thm}

\begin{proof}
The first group of statements are straightforward via appropriate
row and column operations on the Hankel matrix.

For the limits, consider
\[
\lim_{\substack{\beta\to 1\\ |\beta|<1}}
(1-\beta^2) \Exp_{U\in O^+(2l+2)} \det(1-\beta U)^{-1} \det(g(U))
\]
Here we need to compute the limit of
\[
\beta^{-l}
\langle 2^l p^{--}_l(x),(1-\beta^2)/(1-2\beta x+\beta^2)\rangle.
\]
Now,
\[
(1-\beta^2)/(1-\beta (z+1/z)+\beta^2)
=
1+
\sum_{m>0} \beta^m (z^m+z^{-m}),
\]
so for any polynomial $q(z)$,
\[
\lim_{\substack{\beta\to 1\\ |\beta|<1}}
\langle q(z),1/(1-2\beta x+\beta^2)\rangle
=
q(1) g(1)^2.
\]
Plugging in and simplifying gives the desired result.  The other cases are
either analogous or follow from Theorem \ref{thm:alphaformulae} above.
\end{proof}

\section{Diagonal points}\label{sec:diags}

Recall that in defining the ensembles $S^\symmg$ above, we insisted
for $\tsymmO$ and $\tsymmu$ that $\pi$ have no fixed points,
and for $\tsymmS$ and $\tsymmu$ that $\pi\iota$ have no fixed points.
While these conditions are very natural from the standpoint of
points in the square, they seem somewhat artificial in the permutation
setting.  This suggests consideration of the following extended ensembles:
{
\allowdisplaybreaks
\begin{align}
\tilde{S}^{\symmO}_n &= \{\pi\in S_n:\pi=\pi^{-1}\}\\
\tilde{S}^{\symmS}_n &= \{\pi\in S_n:\pi=\iota\pi^{-1}\iota\}\\
\tilde{S}^{\symmu}_n
&=
\{\pi\in S_{2n}:\pi=\pi^{-1},\ \pi=\iota\pi^{-1}\iota\}.
\end{align}
}

We immediately encounter a difficulty, however, with the Poisson
generating function, since the probabilities involve division by
the complicated sum
\[
\sum_{0\le m\le \lfloor n/2\rfloor} \binom{n}{2m} \frac{(2m)!}{2^mm!}
\]
for $\tsymmS$ and $\tsymmO$, and a similar sum for $\tsymmu$.  However,
there is a generalized Poisson model for which the generating function
is again tractable.  For $\tsymmO$, the model is as follows:
In a given infinitesimal time interval
$[t,t+dt]$, we add two points (the images of a uniform random point) with
probability $(1/2)tdt$ and add one point (a point uniformly chosen on 
the $\diagup$ diagonal) with probability $\alpha dt$.  Then the probability
that we have $n$ points at time $t$ is
\[
e^{-\alpha t-t^2/2}
\left(\sum_{0\le m\le \lfloor n/2\rfloor} \alpha^{n-2m} \binom{n}{2m} \frac{(2m)!}{2^m m!}\right)
t^n/n!.
\]
Thus if we define $P^\symmO_l(t;\alpha)$ to be the probability
that the increasing subsequence length at time $t$ is at most $l$, then
\[\label{diags3.6}
P^\symmO_l(t;\alpha)
=
e^{-\alpha t-t^2/2}
\sum_{0\le n} \frac{t^n}{n!}
\sum_{0\le m}
\alpha^m \tilde{f}^\symmO_{nml}
,
\]
where $\tilde{f}^\symmO_{nml}$ is the number of involutions on $n$ elements
with $m$ fixed points and no increasing subsequence of length greater than
$l$.  This is, of course, compatible with our previous notation, with
$P^\symmO_l(t;0)=P^\symmO_l(t)$.
Similarly, for $\tsymmS$, if we add negated points with probability
$\beta dt$, we obtain the Poisson generating function
\[\label{diags3.7}
P^\symmS_l(t;\beta)
=
e^{-\beta t-t^2/2}
\sum_{0\le n} \frac{t^n}{n!}
\sum_{0\le m}
\beta^m \tilde{f}^\symmS_{nml}.
\]
Finally, for $\tsymmu$, we have two parameters in the model; we add
fixed points with probability $\alpha dt$, negated points with
probability $\beta dt$, and generic points with probability $2tdt$,
obtaining
\[\label{diags3.8}
P^\symmu_l(t;\alpha,\beta)
=
e^{-\alpha t-\beta t-t^2}
\sum_{0\le n} \frac{t^n}{n!}
\sum_{0\le m_+,m_-}
\alpha^{m_+} \beta^{m_-} \tilde{f}^\symmu_{nm_+m_-l}.
\]
The $\beta$ parameter is largely irrelevant, since:

\begin{lem}
For all $l\ge 0$, $\alpha\ge 0$ and $\beta\ge 0$,
{
\begin{align}
P^\symmS_{2l+1}(t;\beta)&=P^\symmS_{2l}(t;0)\\
P^\symmu_{2l+1}(t;\alpha,\beta)&=P^\symmu_{2l}(t;\alpha,0)
\end{align}
}
\end{lem}

\begin{proof}
The key observation is that if an increasing subsequence contains a given
point off the $\diagdown$ diagonal, it can always be extended to contain
the reflection of that point through that diagonal; moreover, no increasing
subsequence can contain more than one point on that diagonal.  Thus the
point collection at time $t$ has increasing subset size at most $2l+1$ if
and only if the off-diagonal subset has increasing subset size at most
$2l$.  Since the points were added via a generalized Poisson process, the
off-diagonal subset itself corresponds to a generalized Poisson process;
the lemma is then immediate.
\end{proof}

We can again express these generating functions as integrals:

\begin{thm}\label{thm:fixedpoissons}
For all $l$,
\begin{align}
P^\symmO_l(t;\alpha) &=
e^{-\alpha t-t^2/2} \Exp_{U\in O(l)} \det(1+\alpha U) \exp(t\Tr(U))\\
P^\symmO_l(t;1) &=
e^{-t^2/2} \Exp_{U\in O^-(l+1)} \exp(t\Tr(U))\\
P^\symmS_{2l}(t;\beta) &=
e^{-\beta t-t^2/2} \Exp_{U\in Sp(2l)} \det(1-\beta U)^{-1} \exp(t\Tr(U))
\label{eq:badbeta1}\\
P^\symmS_{2l+1}(t;\beta) &=
e^{-t^2/2} \Exp_{U\in Sp(2l)} \exp(t\Tr(U))\\
P^\symmS_l(t;1) &=
e^{-t^2/2} \Exp_{U\in O^-(l+1)} \exp(t\Tr(U))\\
P^\symmu_{2l}(t;\alpha,\beta) &=
e^{-\alpha t-\beta t-t^2} \Exp_{U\in U(l)} \det((1+\alpha U)(1-\beta
U)^{-1}) \exp(2t\Re\Tr(U))\label{eq:badbeta2}\\
P^\symmu_{2l+1}(t;\alpha,\beta) &=
e^{-\alpha t-t^2} \Exp_{U\in U(l)} \det(1+\alpha U) \exp(2t\Re\Tr(U))\\
P^\symmu_{2l+1}(t;1,\beta) &=
e^{-t^2} \Exp_{U\in O^-(l+2)} \exp(t\Tr(U))
\cdot\Exp_{U\in O^-(l+1)} \exp(t\Tr(U)).
\end{align}
\end{thm}

\begin{proof}
This follows immediately from the symmetric function identities of
Section~\ref{sec:symmfunc} below.  (See, in particular, remark
\ref{rem:tforms} to Corollary \ref{cor:symmfuncdets1}.)
\end{proof}

\begin{rem}
Strictly speaking, equations \eqref{eq:badbeta1} and \eqref{eq:badbeta2}
are only valid when $\beta<1$.  The correct formulae for $\beta\ge 1$
are obtained by analytic continuation (e.g., by expanding $\det(1-\beta
U)^{-1}$ in a formal power series, integrating term-by-term, and summing).
\end{rem}

From Theorem \ref{thm:fixedpoissons} and the orthogonal polynomial
identities of Section~\ref{sec:orthopoly} below, we have the following
formulae.

\begin{cor}\label{cor12}
For $\alpha,\beta \geq 0$, 
\begin{align}
   P_{2l}^\symmO(t;\alpha) &= e^{-\alpha t-t^2/2}
{1\over 2}
\biggl\{ \bigl[\pi^*_{2l-1}(-\alpha;t)
-\alpha\pi_{2l-1}(-\alpha;t)\bigr] D_l^{--}(t)
+ \bigl[\pi^*_{2l-1}(-\alpha;t)
+\alpha\pi_{2l-1}(-\alpha;t)\bigr] D_{l-1}^{++}(t) 
\biggr\}, \\
   P_{2l+1}^\symmO(t;\alpha) &= e^{-\alpha t-t^2/2}
{1\over 2}
\biggl\{ \bigl[\pi^*_{2l}(-\alpha;t)+\alpha\pi_{2l}(-\alpha;t)\bigr] 
e^{t}D_l^{+-}(t)
+  \bigl[\pi^*_{2l}(-\alpha;t)-\alpha\pi_{2l}(-\alpha;t)\bigr] 
e^{-t}D_{l}^{-+}(t) \biggr\},\\
   P_{2l}^\symmO(t;1) &= P_{2l}^\symmS(t;1)
= e^{-t-t^2/2} D_l^{-+}(t),\\
   P_{2l+1}^\symmO(t;1) &= P_{2l+1}^\symmS(t;1)
= e^{-t^2/2} D_l^{++}(t),\\
   P_{2l+1}^\symmS(t;\beta) &= e^{-t^2/2} D_l^{++}(t),\\
   P_{2l+1}^\symmu(t;\alpha,\beta) &= e^{-\alpha t-t^2}
\pi^*_l(-\alpha;t) D_l(t),\\
   P_{4l+1}^\symmu(t;1,\beta) &= e^{-t-t^2}
D_l^{++}(t)D_l^{-+}(t),\\
   P_{4l+3}^\symmu(t;1,\beta) &= e^{-t-t^2}
D_l^{++}(t)D_{l+1}^{-+}(t) .
\end{align}
\end{cor}

\section{Schur function identities}\label{sec:symmfunc}

In order to derive the above generating functions for $P^\symmO(t;\alpha)$,
$P^\symmS(t;\beta)$ and $P^\symmu(t;\alpha,\beta)$, we will first
need a stronger version of Theorem~\ref{thm:trmoments}.

\begin{rem} In the sequel, we will see a number of expressions of the
form
\[
\Exp_{U\in U^\symmg(l)} \det(g(U)),
\]
where $g$ takes values in a ring of formal power series (e.g., the ring of
symmetric functions) with coefficients in ${\mathbb C}[z,1/z]$ (Laurent
polynomials in $z$).  This {\it formal} integral is defined by expanding
$\det(g(U))$ and integrating term by term.  We can recover analytical
results by specializing down to a small number of variables in such
a way that the resulting series converge.
\end{rem}

We refer the reader to \cite{Macdonald} for an introduction to symmetric
functions.  Recall the following symmetric function identities:

\begin{lem}
The following identities hold in the algebra of symmetric functions on
two sets of variables:
{
\begin{align}
\prod_{j,k} (1-x_j y_k)^{-1}
&=
\sum_\lambda s_\lambda(x)s_\lambda(y)\\
&=
\exp\left(\sum_{m>0} {1\over m} p_m(x) p_m(y)\right).
\end{align}
}
And for any symmetric function $f$,
\[
\langle f(x),\prod_{j,k} (1-x_j y_k)^{-1} \rangle
=
f(y).
\]
\end{lem}

It will be convenient to write
\[
\prod_{j,k} (1-x_jy_k)^{-1} = \prod_j H(x_j;y),
\]
where
\[
H(u;y) \defeq \prod_j (1-u y_j)^{-1} = \sum_j h_j u^j
\]
is the generating function for the complete symmetric functions $h_j$;
by convention, $h_j=0$ when $j<0$.  In the sequel, we will also
need the generating function
\[
E(u;y) \defeq \prod_j (1+u y_j) = \sum_j e_j u^j
\]
for the elementary symmetric functions.

\begin{thm}
For all $l\ge 0$,
\begin{align}
\sum_{\ell(\lambda)\le l} s_\lambda(x) s_\lambda(y)
&=
\Exp_{U\in U(l)} \det(H(U;x)H(U^\dagger;y))\\
\sum_{\ell(\lambda)\le l} s_{2\lambda}(x)
&=
\Exp_{U\in O(l)} \det(H(U;x))\\
\sum_{\ell(\lambda)\le 2l} s_{\lambda^2}(x)
&=
\Exp_{U\in Sp(2l)} \det(H(U;x))
\end{align}
\end{thm}

\begin{proof}
Consider the first identity.  We have
\[
\Exp_{U\in U(l)} \det(H(U;x)H(U^\dagger;y))
=
\sum_{\lambda,\mu}
s_\lambda(x) s_\mu(y)
\Exp_{U\in U(l)} s_\lambda(U) s_\mu(U^\dagger).
\]
But the integral on the right is 1 when $\lambda=\mu$ and
$\ell(\lambda),\lambda(\mu)\le l$, and is 0 otherwise.  The result follows
immediately.

Similarly, the identity for $O(l)$ follows from the fact that
\[
\Exp_{U\in O(l)} s_\lambda(U)
\]
is 1 when $\lambda$ is even and $\ell(\lambda)\le l$, and is 0 otherwise;
and the identity for $Sp(2l)$ follows from the fact that
\[
\Exp_{U\in Sp(2l)} s_\lambda(U)
\]
is 1 when $\lambda'$ is even and $\ell(\lambda)\le l$ and is 0 otherwise.
\end{proof}

\begin{rems} This is a direct generalization of the argument in
\cite{Rains}.
\end{rems}
\begin{rems} The expression $\det(H(U;x))$ can also be written as
\[
\exp(\sum_{m>0} {1\over m} p_m(x) \Tr(U^m)).
\]
\end{rems}

For the hyperoctahedral cases $\tsymmUU$ and $\tsymmu$, we will need the
notation
\[
\tilde{s}_\lambda(x) = s_{\lambda^{(0)}}(x) s_{\lambda^{(1)}}(x)
\]
whenever $\lambda$ is a partition with trivial 2-core and 2-quotient
$(\lambda^{(0)},\lambda^{(1)})$, and $\tilde{s}_\lambda(x)=0$
otherwise.  From Ex. 1.5.24 of \cite{Macdonald}, this can also be written
as
\[
\tilde{s}_\lambda(x) = (-1)^{f(\lambda)/2} \phi_2(s_\lambda(x)),
\]
where $f(\lambda)$ is the number of odd parts of $\lambda$ (note that this
is even when $\lambda$ has trivial $2$-core), and $\phi_2$ is the
homomorphism such that
\begin{align}
\phi_2(h_{2n}) &= h_n\\
\phi_2(h_{2n+1}) &= 0,
\end{align}
or equivalently,
\[
\phi_2(H(u;x)) = H(u^2;x).
\]

\begin{cor}
For all $l\ge 0$,
\begin{align}
\sum_{\ell(\lambda)\le 2l}
\tilde{s}_\lambda(x)\tilde{s}_\lambda(y)
&=
\Exp_{U\in U^\symmUU(2l)} \det(H(U;x)H(U^\dagger;y))\\
\sum_{\ell(\lambda)\le l} \tilde{s}_{2\lambda^2}(x)
&=
\Exp_{U\in U^\symmu(2l)} \det(H(U;x))
\end{align}
\end{cor}

\begin{proof}
This follows from the computations
\begin{align}
\sum_{\ell(\lambda)\le 2l}
\tilde{s}_\lambda(x)\tilde{s}_\lambda(y)
&=
\sum_{\ell(\lambda)\le l,\ \ell(\mu)\le l}
s_\lambda(x) s_\mu(x) s_\lambda(y) s_\mu(y)\\
&=
\left(\sum_{\ell(\lambda)\le l} s_\lambda(x)s_\lambda(y)\right)^2\\
&=
(\Exp_{U\in U(l)} \det(H(U;x) H(U^\dagger;y)))^2\\
&=
\Exp_{U_1,U_2 \in U(l)} \det(H(U_1;x) H(U_1^\dagger;y))\det(
H(U_2;x) H(U_2^\dagger;y))\\
&=
\Exp_{U\in U(l)\oplus U(l)} \det(H(U;x) H(U^\dagger;y))
\end{align}
and
\begin{align}
\sum_{\ell(\lambda)\le l} \tilde{s}_{2\lambda^2}(x)
&=
\sum_{\ell(\lambda)\le l} s_\lambda(x) s_\lambda(x)\\
&=
\Exp_{U\in U(l)} \det(H(U;x))\det(H(U^\dagger;x))\\
&=
\Exp_{U\in U^\symmu(2l)} \det(H(U;x)).
\end{align}
\end{proof}

The point of writing the integrands in the form $\det(H(U;x))$ is
that we know how to convert such integrals into determinants:

\begin{thm}
For all $l\ge 0$,
\begin{align}
\Exp_{U\in U(l)} \det(H(U;x)H(U^\dagger;y))
&=
\det(g_{j-k}(x;y))_{0\le j,k<l}\\
\Exp_{U\in O^+(2l)} \det(H(U;x))
&=
{1\over 2}
\det(i_{j-k}(x)+i_{j+k}(x))_{0\le j,k<l}\\
\Exp_{U\in O^-(2l)} \det(H(U;x))
&=
H(1;x)H(-1;x)
\det(i_{j-k}(x)-i_{j+k+2}(x))_{0\le j,k<l-1}\\
\Exp_{U\in O^+(2l+1)} \det(H(U;x))
&=
H(1;x)
\det(i_{j-k}(x)-i_{j+k+1}(x))_{0\le j,k<l}\\
\Exp_{U\in O^-(2l+1)} \det(H(U;x))
&=
H(-1;x)
\det(i_{j-k}(x)+i_{j+k+1}(x))_{0\le j,k<l}\\
\Exp_{U\in Sp(2l)} \det(H(U;x))
&=
\det(i_{j-k}(x)-i_{j+k+2}(x))_{0\le j,k<l},
\end{align}
where we define
\begin{align}
g_j(x;y) &= \sum_m h_m(x) h_{m+j}(y)\\
i_j(x) &= \sum_m h_m(x) h_{m+j}(x).
\end{align}
\end{thm}

\begin{proof}
By Theorems \ref{thm2.1} and \ref{thm:intstodets}, it suffices to evalute
\[
{1\over 2\pi}
\int_{[0,2\pi]}
H(e^{i\theta};x)
H(e^{-i\theta};y)
e^{ij\theta}
d\theta.
\]
But this integral is just the coefficient of $z^{-j}$ in
$H(z;x)H(z^{-1};y)$; that is,
\[
\sum_m h_m(x) h_{m+j}(y)=g_j(x;y).
\]
\end{proof}

The following is then immediate:

\begin{cor}\label{cor:symmfuncdets1}
For all $l>0$,
\begin{align}
\sum_{\ell(\lambda)\le l} s_\lambda(x) s_\lambda(y)
&=
\det(g_{j-k}(x;y))_{0\le j,k<l}\\
\sum_{\ell(\lambda)\le 2l} s_{2\lambda}(x)
&=\left[
{1\over 2}\det(i_{j-k}(x)+i_{j+k}(x))_{0\le j,k<l}
+
H(1;x)H(-1;x)
\det(i_{j-k}(x)-i_{j+k+2}(x))_{0\le j,k<l-1}
\right]/2\label{eq:5.41}\\
\sum_{\ell(\lambda)\le 2l+1} s_{2\lambda}(x)
&=\left[
H(1;x)
\det(i_{j-k}(x)-i_{j+k+1}(x))_{0\le j,k<l}
+
H(-1;x)
\det(i_{j-k}(x)+i_{j+k+1}(x))_{0\le j,k<l}
\right]/2\label{eq:5.42}\\
\sum_{\ell(\lambda)\le 2l} s_{\lambda^2}(x)
&=
\det(i_{j-k}(x)-i_{j+k+2}(x))_{0\le j,k<l}
\end{align}
\end{cor}

\begin{rems}
The expressions for $\sum_{\ell(\lambda)\le l} s_\lambda(x)s_\lambda(y)$
and $\sum_{\ell(\lambda)\le 2l} s_{\lambda^2}(x)$ are already known
(\cite{Gessel} and \cite{Goulden}, respectively; the latter also gives
analogues of \eqref{eq:5.64} and \eqref{eq:5.65} below), although it is
worth noting that the proof of the latter formula given here, in contrast
to most of the proofs in the literature, involves neither pfaffians nor
symplectic tableaux.  The particular forms of \eqref{eq:5.41} and
\eqref{eq:5.42} given here are new, although (C. Krattenthaler, personal
communication) they can be readily derived from Theorem 2.4(3) of
\cite{Okada} together with well-known character formulas for symplectic and
special orthogonal characters; see Remark \ref{rem:chars} below.  The
analogue of \eqref{eq:5.63}, again given in terms of orthogonal and
symplectic characters, is given in Theorem 2 of \cite{Krattenthaler2}.
\end{rems}
\begin{rems}\label{rem:tforms}
We recover our earlier identities by noting that when $|\lambda|=n$,
$\langle p_{1^n}, s_\lambda(x)\rangle$
is equal to the number of tableaux of shape $\lambda$.  So the
generating functions for $f^\symmg$ are obtained by taking inner products
of the above formulae with $\exp(t p_1(x))$.
But in general taking an inner product with the function
\[
\exp(\sum_{m>0} {f_m\over m} p_m(x))
\]
gives a homomorphism from the algebra of symmetric functions to
${\mathbb Q}[f_1,f_2,\ldots]$, given by the substitutions $p_m(x)\mapsto f_m$.

In our case, we have:
\begin{align}
h_m(x) &\mapsto {t^m\over m!}\\
H(u;x) &\mapsto e^{ut}\\
i_j(x) &\mapsto I_j(2t)
\end{align}
so, for instance, we obtain
\begin{align}
\sum_{n\ge 0} f_{nl} {t^{2n}\over n!^2}
&=
\Exp_{U\in U(l)} \exp(2t\Re\Tr(U))\\
&=
\det(I_{j-k}(2t))_{0\le j,k<l},
\end{align}
as before.
\end{rems}

\begin{rems}\label{rem:chars}
In \cite{Stembridge}, it is observed that the determinant
\[
\det(i_{j-k}(x)+i_{j-k+2}(x))_{0\le j,k<l}
\]
is related by duality to a certain irreducible character of the symplectic
group, and similarly the determinants of Corollary \ref{cor:sfallinv} below
are related to irreducible characters of the (odd-dimensional) special
orthogonal group, as are those of equations \eqref{eq:5.41} and
\eqref{eq:5.42}. It is not clear how this relates to our integrals, since
the integrals are over groups of fixed dimension, while the characters are
associated to groups of arbitrarily high dimension.
\end{rems}

\begin{rems}
The approach in \cite{Goulden} for deriving such formulae is based
on the operator
\[
f(z) \mapsto \phi\Bigl(\prod_{1\le j<k\le l} (1-z_j^{-1} z_k) f(z)\Bigr)
\]
where $\phi$ is the linear operator on Laurent series taking $\prod_{1\le
j\le l} z_j^{b_j}$ to $\prod_{1\le j\le l} h_{b_j}(x)$.  Now, this operator
can be expressed as an integral:
\[
\phi\Bigl(\prod_{1\le j<k\le l}(1-z_j^{-1} z_k) f(z)\Bigr)
=
\int_{[0,2\pi]^l} f(e^{i\theta})
\prod_{1\le j\le l} H(e^{-i\theta_j};x)
\prod_{1\le j<k\le l} (1-e^{-i\theta_j}e^{i\theta_k})
\prod_{1\le j\le l} d\theta_j,
\]
which after symmetrizing the integrand becomes simply
\[
\Exp_{U\in U(l)} f(U) \det(H(U^\dagger;x)).
\]
\end{rems}

In \cite{Goulden}, formulae are given for a generalization of
$\sum_{\ell(\lambda)\le l} s_{\lambda^2}(x)$, in which the sum is over all
$\lambda$ with a given number of odd parts in the dual.  This is derived
using the homomorphism $H^\perp(\beta;x)$ (\cite{Macdonald}) which takes
$f(x)$ to $f(\beta,x)$ for any symmetric function $f$.  The relevant
property of $H^\perp(\beta;x)$ is that
\[
H^\perp(\beta;x) s_\lambda(x)
=
s_\lambda(\beta,x)
=
\sum_{\mu\lesssim\lambda} \beta^{|\lambda|-|\mu|} s_\mu(x),
\]
where $\mu\lesssim\lambda$ indicates that $\lambda-\mu$ is a ``horizontal
strip''; that is,
\[
\lambda'_i\le \mu'_i \le \lambda'_i+1
\]
for all $i$, or equivalently
\[
\lambda_1\ge \mu_1\ge \lambda_2\ge \mu_2\ge \ldots.
\]
We will also need the dual operator $E^\perp(\alpha;x)$, which
satisfies
\[
E^\perp(\alpha;x) s_\lambda(x)
=
\sum_{\mu\lesssim'\lambda} \alpha^{|\lambda|-|\mu|} s_\mu(x).
\]
where
\[
\mu\lesssim' \lambda \Leftrightarrow \mu'\lesssim\lambda'.
\]
We have the identities
\begin{align}
H^\perp(\beta;x)H(u;x) &= (1-\beta u)^{-1} H(u;x)\\
H^\perp(\beta;x)E(u;x) &= (1+\beta u) E(u;x)\\
E^\perp(\alpha;x)H(u;x) &= (1+\alpha u) H(u;x)\\
E^\perp(\alpha;x)E(u;x) &= (1-\alpha u)^{-1} E(u;x).
\end{align}

\begin{thm}\label{thm:fixinv}
For all $l>0$, we have the integral identities
\begin{align}
\sum_{\ell(\lambda)\le l} \alpha^{f(\lambda)} s_\lambda(x)
&=
\Exp_{U\in O(l)} \det((1+\alpha U)H(U;x))\label{eq:5.63}\\
\sum_{\ell(\lambda)\le 2l} \beta^{f(\lambda')} s_\lambda(x)
&=
\Exp_{U\in Sp(2l)} \det((1-\beta U)^{-1}H(U;x))\label{eq:5.64}\\
\sum_{\ell(\lambda)\le 2l+1} \beta^{f(\lambda')} s_\lambda(x)
&=
H(\beta;x) \Exp_{U\in Sp(2l)} \det(H(U;x))\label{eq:5.65}
\end{align}
\end{thm}

\begin{proof}
The point is that given any partition $\lambda$, there is a unique
partition $\mu$ such that $\lambda\lesssim' 2\mu$; simply add one to each
odd part of $\lambda$, then divide by 2.  In particular, then,
$f(\lambda)=|2\mu|-|\lambda|$, and thus
\begin{align}
\sum_{\ell(\lambda)\le l} \alpha^{f(\lambda)} s_\lambda(x)
&=
\sum_{\ell(\mu)\le l}
\sum_{\lambda\lesssim' 2\mu}
\alpha^{|2\mu|-|\lambda|} s_\lambda(x)\\
&=
E^\perp(\alpha;x)
\sum_{\ell(\mu)\le l} s_{2\mu}(x)\\
&=
\Exp_{U\in O(l)}
E^\perp(\alpha;x)
\det(H(U;x))\\
&=
\Exp_{U\in O(l)}
\det((1+\alpha U)H(U;x)).
\end{align}
Similarly,
\begin{align}
\sum_{\ell(\lambda)\le 2l} \beta^{f(\lambda')} s_\lambda(x)
&=
H^\perp(\beta;x) \Exp_{U\in Sp(2l)} \det(H(U;x))\\
&=
\Exp_{U\in Sp(2l)} \det((1-\beta U)^{-1}H(U;x)).
\end{align}
By the argument in \cite{Goulden}, the expression for
\[
\sum_{\ell(\lambda)\le 2l+1} \beta^{f(\lambda')} s_\lambda(x)
\]
follows from Pieri's formula
\[
H(\beta;x) s_\mu(x) = \sum_{\lambda\gtrsim \mu} \beta^{|\lambda|-|\mu|}
s_\lambda(x)
\]
(eq. (5.16) in Chapter I of \cite{Macdonald}), 
and the fact that for every $\lambda$
with $\ell(\lambda)\le 2l+1$, there is a unique $\mu$ with
$\mu^2\lesssim\lambda$, and for that $\mu$, $\ell(\mu)\le l$.
\end{proof}

\begin{rem}
It is worth stressing here that the integral
\[
\Exp_{U\in Sp(2l)} \det((1-\beta U)^{-1}H(U;x))
\]
is a {\it formal} integral, in which the expression
$
(1-\beta z)^{-1}
$
stands for the formal power series
$
\sum_{0\le m} z^m \beta^m.
$
We must therefore take special care when specializing to an analytical
integral with $|\beta|\ge 1$.
\end{rem}

\begin{cor}\label{cor:sfallinv}
For all $l>0$,
\[
\sum_{\ell(\lambda)\le l} s_\lambda(x)
=
E(1;x) \Exp_{U\in O^-(l+1)} \det(H(U;x)),
\]
so
\begin{align}
\sum_{\ell(\lambda)\le 2l} s_\lambda(x)
&=
\det(i_{j-k}(x)+i_{j+k+1}(x))_{0\le j,k<l}\\
\sum_{\ell(\lambda)\le 2l+1} s_\lambda(x)
&=
H(1;x)
\det(i_{j-k}(x)-i_{j+k+2}(x))_{0\le j,k<l}
\end{align}
\end{cor}

\begin{proof}
We have
\begin{align}
\sum_{\ell(\lambda)\le l} s_\lambda(x)
&=
\sum_{\ell(\lambda)\le l} 1^{f(\lambda)} s_\lambda(x)\\
&=
\Exp_{U\in O(l)} \det((1+U)H(U;x))\\
&=
{1\over 2} \Exp_{U\in O^+(l)} \det((1+U)H(U;x))
\end{align}
since $\det(1+U)$ vanishes on $O^-(l)$.  But then we can apply Theorem
\ref{thm:alphaformulae} to obtain
\[
{1\over 2} \Exp_{U\in O^+(l)} \det((1+U)H(U;x))
=
H(-1;x)^{-1} \Exp_{U\in O^-(l+1)} \det(H(U;x))
\]
as desired (recall $H(-t;x)E(t;x)=1$).  The remaining formulas are immediate.
\end{proof}

\begin{rems}
Again, the determinantal forms are known \cite{Gordon:NOPPV} (based
on an earlier pfaffian form \cite{GordonHouten}).
\end{rems}
\begin{rems}
We could also
have derived the formulae by the (formal) substitution $\beta=1$
in the above symplectic formulae, or by taking the limit $\beta\to 1^{-}$.
\end{rems}

We also have the following formulae:

\begin{thm}\label{thm:sfeig2inv}
For all $l\ge 0$,
\begin{align}
\sum_{\lambda'_2\le l} \alpha^{f(\lambda)} s_\lambda(x)
&=
E(\alpha;x)E_{U\in O(l)} \det(H(U;x)) 
\label{eq:sfeig2inv1}
\\
\sum_{\lambda'_2\le l\le \lambda'_1} \alpha^{2\lambda'_1-l-f(\lambda)} s_\lambda(x)
&=
E(\alpha;x)E_{U\in O(l)} \det(U)\det(H(U;x)).
\label{eq:sfeig2inv2}
\end{align}
\end{thm}

\begin{proof}
We recall
\[
E(\alpha;x) s_\mu(x) = \sum_{\mu\lesssim'\lambda} \alpha^{|\lambda|-|\mu|}
s_\lambda(x);
\]
that is, we sum over all partitions $\lambda$ obtained from $\mu$ by
adjoining a vertical strip.  On the other hand,
\begin{align}
E_{U\in O(l)} \det(H(U;x))
&=
\sum_{\substack{\ell(\mu)\le l\\f(\mu)=0}}
s_\mu(x)\\
E_{U\in O(l)} \det(U) \det(H(U;x))
&=
\sum_{\ell(\mu)=f(\mu)=l}
s_\mu(x).
\end{align}

In each case, there is at most one way to remove a vertical strip from a
generic partition to obtain one with the desired special form.  Thus it
remains only to determine which partitions occur.  In either case, we note
that $\mu_{l+1}=0$ so $\lambda_{l+1}\le 1$; it follows that $\lambda'_2\le l$.
In the second case, $\mu_l\ge 1$, and thus $\lambda'_1\ge l$.  As these
conditions are readily seen to be sufficient, the result follows.
\end{proof}

\begin{rem}
This gives another proof of Corollary \ref{cor:sfallinv}, by 
subtracting \eqref{eq:sfeig2inv2} from \eqref{eq:sfeig2inv1} with $\alpha=1$.
\end{rem}

It remains to consider the formulae corresponding to hyperoctahedral
involutions.  For an arbitrary partition $\lambda$, we define new partitions
\begin{align}
\lambda^+ &= (\lambda_1,\lambda_3,\ldots)\\
\lambda^- &= (\lambda_2,\lambda_4,\ldots).
\end{align}
Note that $\lambda^-$ and $\lambda^+$ are the unique partitions such that
$(\lambda^-)^2\lesssim\lambda\lesssim(\lambda^+)^2$; also note that
$f(\lambda')=|\lambda^+|-|\lambda^-|$.

\begin{lem}
A partition $\lambda$ has trivial 2-core if and only if
$f(\lambda^+)=f(\lambda^-)=f(\lambda)/2$.
\end{lem}

\begin{proof}
A partition has trivial 2-core if and only if its diagram can be tiled by
dominos.  By a classical result, this can happen if and only if the diagram
contains as many points $(i,j)$ with $i+j$ even as with $i+j$ odd.
For $a$, $b\in \{0,1\}$, let $C_{ab}$ be the number of points in the
diagram of $\lambda$ with $(i\bmod 2,j\bmod 2)=(a,b)$.  Then
\begin{align}
f(\lambda^+)&=C_{11}-C_{10}\\
f(\lambda^-)&=C_{01}-C_{00}.
\end{align}
But then
\[
f(\lambda^+)-f(\lambda^-) = C_{11}-C_{10}-C_{01}+C_{00} = 0.
\]
\end{proof}

\begin{thm}
For all $l\ge 0$,
\[
\sum_{\ell(\lambda)\le 2l} 
\alpha^{f(\lambda)/2} \beta^{f(\lambda')/2}
\tilde{s}_\lambda(x)
=
\Exp_{U\in U(l)} \det((1+\alpha U)(1-\beta U)^{-1}H(U;x)H(U^\dagger;x)).
\]
\end{thm}

\begin{proof}
It follows from the lemma that for any $\lambda$ with
trivial 2-core, there exist unique partitions $\mu$ and $\nu$ such that
\[
\lambda\lesssim \nu^2 \lesssim' 2\mu^2,
\]
and which satisfy
\begin{gather}
f(\lambda) = 2 f(\nu) = |2\mu^2|-|\nu^2|\\
f(\lambda')= |\nu^2|-|\lambda|
\end{gather}
Consequently,
\[
\phi_2\left(
\sum_{\ell(\lambda)\le 2l} 
(-\alpha)^{f(\lambda)/2}
\beta^{f(\lambda')/2}
s_\lambda\right)
=
\phi_2\left(
\sum_{\ell(\mu)\le l}
E^\perp(\sqrt{-\alpha};x)H^\perp(\sqrt{\beta};x) s_{2\mu^2}(x)\right).
\]
It thus suffices to show that for any partition $\mu$,
\[
\phi_2(E^\perp(a;x)H^\perp(b;x) s_{2\mu^2}(x))
=
s_\mu(x) (E^\perp(-a^2;x) H^\perp(b^2;x) s_\mu(x)),
\]
since then we may set $y=x$ in
\[
\sum_{\ell(\mu)\le l}
s_\mu(y) (E^\perp(\alpha;x) H^\perp(\beta;x) s_\mu(x))
=
\Exp_{U\in U(l)} \det((1+\alpha U)(1-\beta U)^{-1}H(U;x)H(U^\dagger;y)).
\]

By the Jacobi-Trudi identity, $s_{2\mu^2}(x)$ is the determinant
of the block matrix
\[
\det(B_{\mu_i+j-i})_{1\le i,j\le \ell(\mu)},
\]
where
\[
B_m = \pmatrix
h_{2m} & h_{2m-1}\\
h_{2m+1}& h_{2m}
\endpmatrix
\]
Equivalently, via row and column operations, we could take
\[
B_m = \pmatrix
h_{2m}-b h_{2m-1} & h_{2m-1}\\
h_{2m+1}-(a+b)h_{2m}+ab h_{2m-1} &
h_{2m}-a h_{2m-1}
\endpmatrix.
\]
Upon applying $E^\perp(a;x)H^\perp(b;x)$ then $\phi_2$,
we obtain
\[
B'_m = \pmatrix
h_m&(a+b) H^\perp(b^2;x) h_{m-1}\\
0&E^\perp(-a^2;x) H^\perp(b^2;x) h_m
\endpmatrix.
\]
But then
\[
\det(B'_{\mu_i+j-i})_{1\le i,j\le \ell(\mu)}
=
s_\mu(x) (E^\perp(-a^2;x) H^\perp(b^2;x) s_\mu(x)),
\]
as required.
\end{proof}

\begin{cor}\label{cor:sfhyper}
\[
\sum_{\ell(\lambda)\le 2l+1} \alpha^{f(\lambda)/2} \beta^{f(\lambda')/2}
\tilde{s}_\lambda(x)
=
H(\beta;x)
\Exp_{U\in U(l)} \det((1+\alpha U)H(U;x)H(U^\dagger;x))
\]
\end{cor}

\begin{proof}
Here we use the following analogue of Pieri's formula:
\begin{align}
H(\beta;x) \tilde{s}_\lambda(x)
&=
(-1)^{f(\lambda)/2}
\phi_2(H(\sqrt{\beta};x) s_\lambda(x))\\
&=
\sum_{\mu\gtrsim\lambda}
(-1)^{(f(\lambda)-f(\mu))/2}
\beta^{(|\lambda|-|\mu|)/2} \tilde{s}_\mu(x).
\end{align}

Thus
\begin{align}
H(\beta;x)
\sum_{\ell(\lambda)\le l} \alpha^{f(\lambda)} \tilde{s}_{\lambda^2}(x)
&=
\sum_{\substack{\ell(\lambda)\le l\\ \mu\gtrsim \lambda^2}}
\beta^{|\mu|/2-|\lambda|}
\alpha^{f(\lambda)} \tilde{s}_{\mu}(x)\\
&=
\sum_{\ell(\mu)\le 2l+1}
\beta^{|\mu|/2-|\mu^-|}
\alpha^{f(\mu^-)} \tilde{s}_{\mu}(x).\\
&=
\sum_{\ell(\mu)\le 2l+1}
\beta^{f(\mu')/2}
\alpha^{f(\mu)/2} \tilde{s}_{\mu}(x).
\end{align}
\end{proof}

Analogously,

\begin{cor}
For all $l>0$,
\begin{align}
\sum_{\mu'_2\le 2l}
\alpha^{f(\mu)/2}
\beta^{f(\mu')/2}
\tilde{s}_\mu(x)
&=
E(\alpha;x)
E_{U\in U(l)} \det((1-\beta U)^{-1} H(U;x) H(U^\dagger;x))\\
\sum_{\mu'_2\le 2l+1}
\alpha^{f(\mu)/2}
\beta^{f(\mu')/2}
\tilde{s}_\mu(x)
&=
E(\alpha;x) H(\beta;x)
E_{U\in U(l)} H(U;x) H(U^\dagger;x))
\end{align}
\end{cor}

\begin{proof}
Dualizing the proof of Corollary \ref{cor:sfhyper}, we find
\[
E(\alpha;x)
\sum_{\ell(\lambda)\le l} \beta^{f(\lambda')} \tilde{s}_{2\lambda}(x)
=
\sum_\mu
\alpha^{f(\mu)/2}
\beta^{f(\mu')/2}
\tilde{s}_\mu(x),
\]
where the sum is over $\mu$ such that $((\mu')^{-})'=\lambda$ with
$\ell(\lambda)\le l$.  But, as in the proof of Theorem \ref{thm:sfeig2inv},
this condition is simply that $\mu'_2\le l$.  The result follows.
\end{proof}

\section{Pfaffian identities}\label{sec:pfaffians}

As we remarked earlier, many of the earlier proofs of the known identities
from Section \ref{sec:symmfunc} proceeded via an intermediate pfaffian
form.  We sketch analogous proofs of the remaining identities of that
section.

We recall the definition of the pfaffian.  If $A$ is a $2n\times 2n$
antisymmetric matrix, the pfaffian $\pf(A)$ is defined by:
\[
\pf(A) = {1\over 2^n n!} \sum_{\pi\in S_{2n}} \sigma(\pi)
\prod_{1\le j\le n} A_{\pi(2j-1)\pi(2j)}.
\]
In the odd-dimensional case, it will be convenient to use the notation
\[
\pf(v;A)
\]
to denote the pfaffian of $A$ bordered by $v$.

The fundamental identity we use is the following:

\begin{thm} (de Bruijn, \cite{deBruijn}) \label{thm:deBruijn}
Let $X$ be a measure space, let $\rho(x,y)$ be an antisymmetric function on
$X\times X$, and let $\phi_j(x)$ be a sequence of functions on $X$,
such that for all $j$ and $k$, the function $\phi_j(x) \rho(x,y) \phi_k(y)$
is absolutely integrable.  Then for $n$ even,
\[
\int_{\vec{x}}
\pf(\rho(x_j,x_k))_{1\le j,k\le n}
\det(\phi_j(x_k))_{1\le j,k\le n}
=
n!
\pf(
\int_{x,y} \phi_j(x) \rho(x,y) \phi_k(y)
)_{1\le j,k\le n}.
\]
For $n$ odd,
\[
\int_{\vec{x}}
\pf(\rho(x_j,x_k))_{1\le j,k\le n}
\det(\phi_j(x_k))_{1\le j,k\le n}
=
n!
\pf(
(\int_x \phi_j(x))_{1\le j\le n};
(\int_{x,y} \phi_j(x) \rho(x,y) \phi_k(y))_{1\le j,k\le n}
),
\]
where the integrals are all over $X$.
\end{thm}

\begin{proof}
For $n$ even, we have:
\begin{align}
\int_{\vec{x}}
\pf(\rho(x_j,x_k))_{1\le j,k\le n}
\det(\phi_j(x_k))_{1\le j,k\le n}
&=
\sum_{\pi_1,\pi_2\in S_n}
\sigma(\pi_1\pi_2)
\int_{\vec{x}}
\prod_{1\le j\le n/2} \rho(x_{\pi_1(2j-1)},x_{\pi_2(2j)})
\prod_{1\le j\le n} \phi_j(x_{\pi_2(j)})\\
&=
n!
\sum_{\pi\in S_n}
\sigma(\pi)
\int_{\vec{x}}
\prod_{1\le j\le n/2} \rho(x_{\pi(2j-1)},x_{\pi(2j)})
\prod_{1\le j\le n} \phi_{\pi(j)}(x_{\pi(j)})\\
&=
n!
\sum_{\pi\in S_n}
\sigma(\pi)
\prod_{1\le j\le n/2}
\int_{x,y} \phi_{\pi(2j-1)}(x)\rho(x,y) \phi_{\pi(2j)}(y).
\end{align}

For $n$ odd, we simply adjoin an extra element $\infty$ to $X$ and
add a new function $\phi_0$ which is 0 away from $\infty$, then
apply the identity for $n$ even.
\end{proof}

\begin{rem}
As noted in \cite{deBruijn}, since the proof holds for arbitrary measure
spaces, this includes the case of a sum, thus obtaining a pfaffian analogue
of the Cauchy-Binet formula (which is also readily obtained as a special
case).  As such, this was independently rediscovered in
\cite{IW}, as well as some extensions (e.g., to a $q$-analogue) which we
will not need.  It is also worth noting that de Bruijn's formula has
applications to random matrices (\cite{TracyWidomcluster}).
\end{rem}

Associated to any partition $\lambda$ of length $\le l$, we associate
a $l$-tuple $\mu_l(\lambda)$ of distinct integers, defined by
\[
\mu_l(\lambda)_j = \lambda_{l-j} +j,
\]
with $0\le j<l$.  In terms of this, the Jacobi-Trudi identity becomes
\[
s_\lambda(x) = \det(h_{\mu_l(\lambda)_k-j}(x))_{0\le j,k<l}.
\]
We observe that this has the appropriate form for the application of
Theorem \ref{thm:deBruijn}, where we take
\[
\phi_j(\mu_k) = \phi_j(\mu_k;x) \defeq h_{\mu_k-j}(x).
\]
Thus, for instance, we have
\begin{align}
\sum_{\ell(\lambda)\le l}
s_\lambda(x) s_\lambda(y)
&=
{1\over l!}
\sum_{\mu_1,\mu_2,\ldots \mu_l\in \N}
\det(\phi_j(\mu_k;x))\det(\phi_j(\mu_k;y))\\
&=
\det(
\sum_{\mu\in \N}
\phi_j(\mu;x)\phi_k(\mu;x)
)_{0\le j,k<l}\\
&=
\det(
g_{j-k}(x;y)
)_{0\le j,k<l}.
\end{align}
This is, of course, precisely Gessel's original proof (\cite{Gessel}),
slightly restated.  Similarly, the identity for $\tsymmUU$ follows
analogously.

To handle the involution cases, we note the following: Define a function
$F(\mu)$ on $l$-tuples of nonnegative integers as follows: if $\mu$ is
increasing (and thus $\mu=\mu_l(\lambda)$ for some $\lambda$), then
\[
F(\mu) = \alpha^{f(\lambda)} \beta^{f(\lambda')}.
\]
Otherwise, $F(\mu)$ is alternating under permutations of $\mu$.

\begin{lem} \cite{IOW}
If $l$ is even, then
\[
F(\mu) = \pf(F(\mu_i\mu_k))_{0\le i,k<l},
\]
while if $l$ is odd, then
\[
F(\mu) = -\pf(F(\mu_i)_{0\le i<l};F(\mu_i\mu_k)_{0\le i,k<l}).
\]
\end{lem}

At this point, we can write
\[
\sum_{\ell(\lambda)\le l}
\alpha^{f(\lambda)} \beta^{f(\lambda')}
s_\lambda(x)
=
{1\over l!}
\sum_{\mu_1,\mu_2,\ldots \mu_l\in \N}
F(\mu)
\det(\phi_j(\mu_k;x))_{0\le j,k<l}.
\]
For $l$ even, this becomes:
\[
{1\over l!}
\sum_{\mu_1,\mu_2,\ldots \mu_l\in \N}
\pf(F(\mu_i\mu_k))_{0\le i,k<l}
\det(\phi_j(\mu_k;x))_{0\le j,k<l}
=
\pf(
\sum_{\mu,\nu} F(\mu\nu) \phi_j(\mu)\phi_k(\nu)
)_{0\le j,k<l}.
\]
For $l$ odd, we must border the pfaffian with
\[
(-\sum_\mu F(\mu) \phi_j(\mu))_{0\le j<l}
\]

In three cases (corresponding to $\tsymmO$, $\tsymmS$, $\tsymmu$), this
pfaffian simplifies.  We consider only $\tsymmO$ in detail (in particular
since the resulting pfaffians for the other cases are more complicated).

For $l$ even, after simplifying, we find
\[
M(j,k)\defeq \sum_{\mu,\nu} F(\mu,\nu) \phi_j(\mu)\phi_k(\nu)
=
M_0(j,k)+M_1(j,k),
\]
where
\[
M_0(j,k)=
{1\over 2}
\sum_{0\le d<(k-j)} 
\left(
(1+\alpha^2) i_{2d+1-(k-j)}(x)
+
\alpha (i_{2d-(k-j)}(x)+i_{2d+2-(k-j)}(x))\right),
\]
while $M_1(j,k)=0$ unless $k-j$ is odd, when
\[
M_1(j,k)=(-1)^j {1-\alpha^2\over 2} H(1;x)H(-1;x)=
(-1)^j H(1;\alpha,x)H(-1;\alpha,x)
\]
Now, $M_1$ has rank 2, so $\pf(M)$ can be simplified via appropriate
row and columns.  To be precise, if we subtract row/column $i$ from
row/column $i+2$ in decreasing order, $M_1$ is now 0 except when
$(j,k)=(0,1)$ or $(1,0)$.  Thus
\[
\sum_{\ell(\lambda)\le l} \alpha^{f(\lambda)} s_\lambda(x)
=
\pf(M_0)+{1\over 2} H(1;\alpha,x)H(-1;\alpha,x)\pf(M'_0),
\]
where
\[
M'_0(j,k) = i_{k-j-1}(x)-i_{k-j+1}(x).
\]
Our earlier formula follows from the following identity due to Gordon
\cite{Gordon:NOPPV}:

\begin{thm}
Let $x_j$ be an odd function of $j\in \Z$.  Then
\begin{align}
\pf(x_{k-j})_{0\le j,k<2l}
&=
\det(
\sum_{0\le m\le k} x_{j-k+2m+1}
)_{0\le j,k<l}
\\
\pf((z^{j-l})_{0\le j<2l+1};(x_{k-j})_{0\le k,j<2l+1})
&=
\det(
\sum_{0\le m\le k} (z+{1\over z}) x_{j-k+2m+1}-(x_{j-k+2m}+x_{j-k+2m+2})
)_{0\le j,k<l}
\end{align}
\end{thm}

\begin{rem}
In addition to the proof in \cite{Gordon:NOPPV}, one can also prove this
by showing that the second pfaffian gives an appropriate orthogonal
polynomial, or by applying a special case of de Bruijn's formula (with
$X$ the unit circle).
\end{rem}

The computation for $l$ odd is analogous.

\section{More increasing subsequence problems}\label{sec:moreiseq}

It turns out that the Schur function sums considered above can in many
cases be given direct interpretations in terms of suitably defined
increasing subsequences of random multisets.  Let $\Omega_1$ and $\Omega_2$
be totally ordered sets, and let $W_1\subset \Omega_1$ and $W_2\subset
\Omega_2$ be arbitrarily chosen subsets, with complements $\overline{W_1}$
and $\overline{W_2}$.  It will also be convenient to use the corresponding
indicator functions $\chi_1$ and $\chi_2$.

\begin{defn}
A $(W_1,W_2)$-increasing sequence in $\Omega_1\times \Omega_2$ is
a sequence $(x_i,y_i)$ such that
\[
x_i \le x_{i+1},\ \text{and}\ y_i \le y_{i+1},
\]
for all $i$, and such that
\begin{align}
x_i=x_{i+1} &\implies x_i\in W_1,\\
y_i=y_{i+1} &\implies y_i\in W_2.
\end{align}
\end{defn}

Thus $W_i$ specifies where the $i$th coordinate is allowed to weakly
increase.  In particular, a $(\Omega_1,\Omega_2)$-increasing sequence is
weakly increasing, while a $(\emptyset,\emptyset)$-increasing sequence is
strictly increasing.  We also have a notion of a $(W_1,W_2)$-decreasing
sequence, in which the inequality for $y$ is reversed.

\begin{defn}
A finite multiset $M$ from $\Omega_1\times\Omega_2$ is
$(W_1,W_2)$-compatible if for all $i$ and $j$ such that
$\chi_1(i)\ne \chi_2(j)$, $(i,j)$ occurs at most once in $M$.
\end{defn}

If $M$ is a $(W_1,W_2)$-compatible multiset from $\Omega_1\times \Omega_2$,
define $l_{W_1W_2}(M)$ to be the length of the longest
$(W_1,W_2)$-increasing subsequence of $M$ (with, of course, the condition
that an element may appear in the sequence at most as many times as it
appears in $M$).  We will also require the notation $l^-_{W_1W_2}(M)$
for the length of the longest $(W_1,W_2)$-decreasing subsequence of $M$.

We now define, for each symmetry type, and for certain choices of
$(W_1,W_2)$ for each symmetry type, a random $(W_1,W_2)$-compatible
multiset, and thus a random longest increasing subsequence length.  We
define the multiset by specifying for each $(x,y)$ the distribution of the
multiplicity $w(x,y)$ of $(x,y)$ in $M$; the multiplicities are independent
unless otherwise specified.  Each multiset distribution will also depend on
certain parameters.  For convenience, we use the following notations for
the distributions which appear: $g(q)$ for the geometric distribution with
parameter $q<1$, $b(q)$ for a variable which is 0 with probability
$1/(1+q)$ and 1 with probability $q/(1+q)$, and $g'(\alpha,q)$ for a
variable with
\begin{align}
\Prob(g'(\alpha,q)=k) = {1-q^2\over 1+\alpha q} \alpha^{k\,\bmod\, 2} q^k.
\end{align}
We denote the set of multisets from a set $S$ by $\Mult(S)$.

\begin{itemize}
\item[$\tsymmU$:]
The parameters are subsets $W,W'\subset \Z^+$ and nonnegative sequences
$q_i$ and $q'_i$ such that: $q_iq'_j<1$ whenever $\chi_1(i)=\chi_2(j)$ and
such that
\[
Z^\symmU_{W,W'}(q;q') \defeq
\prod_{\chi(i)=\chi'(j)} (1-q_iq'_j)
\prod_{\chi(i)\ne \chi'(j)} (1+q_iq'_j)^{-1}
\ne 0.
\]
$M\in\Mult(\Z^+\times\Z^+)$ is chosen as follows:
\begin{align}
\chi(i)=\chi'(j)&:\ w(i,j)\sim g(q_iq'_j)\\
\chi(i)\ne \chi'(j)&:\ w(i,j)\sim b(q_iq'_j)
\end{align}
Denote the variable $l_{WW'}(M)$ by $l^\symmU_{WW'}(q;q')$.
\item[$\tsymmUU$:]
The parameters are subsets $W,W'\subset \Z$ symmetric under negation and
nonnegative sequences $q_i$ and $q'_i$ such that $q_iq'_j<1$ whenever
$\chi(i)=\chi'(j)$ and such that
\[
Z^\symmUU_{WW'}(q;q') \defeq
\prod_{\chi(i)=\chi'(j)} (1-q_iq'_j)^2
\prod_{\chi(i)\ne \chi'(j)} (1+q_iq'_j)^{-2}
\ne 0.
\]
$M\in\Mult(\Z\times\Z)$ is chosen as follows:
$w(i,j)=0$ if $i=0$ or $j=0$.  Otherwise:
\begin{align}
\chi(i)=\chi'(j)&:\ w(i,j)=w(-i,-j)\sim g(q_{|i|}q'_{|j|})\\
\chi(i)\ne \chi'(j)&:\ w(i,j)=w(-i,-j)\sim b(q_{|i|}q'_{|j|})
\end{align}
Denote the variable $l_{WW'}(M)$ by $l^\symmUU_{WW'}(q;q')$.
\item[$\tsymmO$:]
The parameters are a subset $W\subset \Z^+$, a real number $\alpha\ge 0$
and a nonnegative sequence $q_i$ such that $q_i<1$ for all $i$, $\alpha
q_i<1$ when $\chi(i)=1$, and
\[
Z^\symmO_W(q;\alpha) \defeq
\prod_{\chi(i)=1} (1-\alpha q_i)
\prod_{\chi(i)=0} (1+\alpha q_i)^{-1}(1-q_i^2)
\prod_{\substack{i<j\\\chi(i)=\chi(j)}} (1-q_iq_j)
\prod_{\substack{i<j\\\chi(i)\ne\chi(j)}} (1+q_iq_j)^{-1}
\ne 0
\]
$M\in\Mult(\Z^+\times\Z^+)$ is chosen as follows:
For $i\ne j$,
\begin{align}
\chi(i)=\chi(j)&: w(i,j)=w(j,i)\sim g(q_iq_j)\\
\chi(i)\ne \chi(j)&: w(i,j)=w(j,i)\sim b(q_iq_j)
\end{align}
For $i=j$,
\begin{align}
\chi(i)=1&:\ w(i,i)\sim g(\alpha q_i)\\
\chi(i)=0&:\ w(i,i)\sim g'(\alpha,q_i).
\end{align}
Denote the variable $l_{WW}(M)$ by $l^\symmO_W(q;\alpha)$.
\item[$\tsymmS$:]
The parameters are a subset $W\subset \Z^+$, a real number $\beta\ge 0$ and
a nonnegative sequence $q_i$ such that $q_i<1$ for all $i$, $\beta q_i<1$
when $\chi(i)=0$, and
\[
Z^\symmS_W(q;\beta) \defeq
\prod_{\chi(i)=0} (1-\beta q_i)
\prod_{\chi(i)=1} (1+\beta q_i)^{-1}(1-q_i^2)
\prod_{\substack{i<j\\\chi(i)=\chi(j)}} (1-q_iq_j)
\prod_{\substack{i<j\\\chi(i)\ne\chi(j)}} (1+q_iq_j)^{-1}
\ne 0
\]
$M\in\Mult(\Z^+\times\Z^-)$ is chosen as follows:
For $i\ne -j$,
\begin{align}
\chi(i)=\chi(-j)&: w(i,-j)=w(j,-i)\sim g(q_iq_{-j})\\
\chi(i)\ne \chi(-j)&: w(i,-j)=w(j,-i)\sim b(q_iq_{-j})
\end{align}
For $i=-j$,
\begin{align}
\chi(i)=0&:\ w(i,-i)\sim g(\beta q_i)\\
\chi(i)=1&:\ w(i,-i)\sim g'(\beta,q_i).
\end{align}
Denote the variable $l_{W,-W}(M)$ by $l^\symmS_W(q;\beta)$.
\item[$\tsymmu$:]
The parameters are a subset $W\in \Z$ symmetric under negation, real
numbers $\alpha\ge 0$ and $\beta\ge 0$, and a nonnegative sequence $q_i$
such that $q_i,\alpha q_i<1$ when $\chi(i)=1$, $q_i,\beta q_i<1$ when
$\chi(i)=0$, and
\[
Z^\symmu_W(q;\alpha,\beta) \defeq
\prod_{\chi(i)=0} (1-\beta q_i)(1+\alpha q_i)^{-1}
\prod_{\chi(i)=1} (1+\beta q_i)^{-1}(1-\alpha q_i)
Z^\symmU_{WW}(q;q)
\ne 0.
\]
$M\in\Mult(\Z\times\Z)$ is chosen as follows: $w(i,j)=0$ if $i=0$ or
$j=0$.  Otherwise, for $|i|\ne |j|$,
\begin{align}
\chi(i)=\chi(j)&:\ w(i,j)=w(-i,-j)=w(j,i)=w(-j,-i)\sim g(q_{|i|}q_{|j|})\\
\chi(i)\ne\chi(j)&:\ w(i,j)=w(-i,-j)=w(j,i)=w(-j,-i)\sim b(q_{|i|}q_{|j|});
\end{align}
for $i=j$,
\begin{align}
\chi(i)=0&:\ w(i,i)=w(-i,-i)\sim g'(\alpha,q_{|i|})\\
\chi(i)=1&:\ w(i,i)=w(-i,-i)\sim g(\alpha q_{|i|});
\end{align}
and for $i=-j$,
\begin{align}
\chi(i)=0&:\ w(i,-i)=w(-i,i)\sim g(\beta q_{|i|})\\
\chi(i)=1&:\ w(i,-i)=w(-i,i)\sim g'(\beta,q_{|i|}).
\end{align}
Denote the variable $l_{WW}(M)$ by $l^\symmu_W(q;\alpha,\beta)$.
\end{itemize}

\begin{rem}
The conditions on the parameters are simply (1) that the various
probability distributions are defined, and (2) that $M$ be finite with
probability 1.  In each case the quantity $Z^\symmg$ gives the probability
that the multiset is empty.
\end{rem}

To relate these variables to our Schur function sums, we will
need the notion of a ``super-Schur'' (also known as ``hook Schur'')
function.  If $x_i$ and $y_i$ are countable sets of variables, we define
$s_\lambda(x/y)$ to be the image of $s_\lambda$ under the homomorphism
\[
H(t;z) \mapsto H(t;x/y) \defeq H(t;x)E(t;y).
\]
(See Example 3.21 of \cite{Macdonald}, but note that we have used a
slightly different sign convention).  In particular, since this is defined
via a homomorphism, all of our identities are valid for such functions;
e.g.,
\[
\sum_{\ell(\lambda)\le l} s_\lambda(x/y) s_\lambda(z/w)
=
\Exp_{U\in U(l)}
\det(H(U;x)E(U;y)H(U^\dagger;z)E(U^\dagger;w)).
\]

Then the relation to our symmetric function identities is as follows:

\begin{thm}\label{thm:kurtjg}
For any valid choices of parameters,
{
\allowdisplaybreaks
\begin{align}
\Prob(l^\symmU_{WW'}(q;q')\le l)
&=
Z^\symmU_{WW'}(q;q')
\sum_{\ell(\lambda)\le l}
s_\lambda(q_{\overline{W}}/q_W)s_\lambda(q'_{\overline{W'}}/q'_{W'})\\
\Prob(l^\symmUU_{WW'}(q;q')\le l)
&=
Z^\symmUU_{WW'}(q;q')
\sum_{\ell(\lambda)\le l}
\tilde{s}_\lambda(q_{\overline{W}}/q_W)
\tilde{s}_\lambda(q'_{\overline{W'}}/q'_{W'})\\
\Prob(l^\symmO_W(q;\alpha)\le l)
&=
Z^\symmO_W(q;\alpha)
\sum_{\ell(\lambda)\le l}
\alpha^{f(\lambda)} s_\lambda(q_{\overline{W}}/q_W)\\
\Prob(l^\symmS_W(q;\beta)\le l)
&=
Z^\symmS_W(q;\beta)
\sum_{\ell(\lambda)\le l}
\beta^{f(\lambda')} s_\lambda(q_{\overline{W}}/q_W)\\
\Prob(l^\symmu_W(q;\alpha,\beta)\le l)
&=
Z^\symmu_W(q;\alpha,\beta)
\sum_{\ell(\lambda)\le l} \alpha^{f(\lambda)/2} \beta^{f(\lambda')/2}
\tilde{s}_\lambda(q_{\overline{W}}/q_W).
\end{align}
}
\end{thm}

\begin{rem}
These processes are generalizations of processes studied by Johansson in
\cite{kurtj:shape}.  In particular, he studies the process
$l^\symmU_{\Z^+\Z^+}$ in the special case
\[
q_i=\cases \sqrt{q}&1\le i\le N\\
           0&i>N\endcases,
\quad
q'_i=\cases \sqrt{q}&1\le i\le M\\
           0&i>M\endcases,
\]
as well as the process $l^\symmO_{\Z^+}$ in the special case
\[
\alpha = 1,\ 
q_i = \cases \sqrt{q}&1\le i\le N\\
              0 & i>N\endcases.
\]
In both cases, he also studies an appropriate limit as $q\to 1$.
\end{rem}

To prove this theorem, we will need a generalization of the Knuth
correspondences \cite{Knuth:correspondence}.  Let $\Omega$ be a totally
ordered set, let $W$ be a chosen subset, and let $\lambda$ be a partition.

\begin{defn}
An $(\Omega,W)$-bitableau $T$ of shape $\lambda$ is a function from the
diagram of $\lambda$ to $\Omega$ which is weakly increasing along each
row and column, such that any element of $W$ appears at most
once in each column, and such that any element of $\overline{W}$ appears
at most once in each row.
\end{defn}

\begin{rem}
This is essentialy the same as the notion of bitableau given in Example
5.23 of \cite{Macdonald}; see also \cite{BereleRegev}.
\end{rem}

We denote the set of such bitableau as $B_\lambda(\Omega,W)$, and observe
(ibid.) that if $x_i$ is a sequence of indeterminates, then for any subset
$W\in Z^+$,
\[
\sum_{T\in B_\lambda(\Z^+,W)} x^T
=
s_\lambda(x_W/x_{\overline{W}}),
\]
where for a bitableau $T$, $x^T$ is the product of $x^{m_i}_i$, where $i$
appears in $T$ $m_i$ times.

\begin{thm}
Given a pair $\Omega_1$ and $\Omega_2$ of totally ordered sets, and
respective subsets $W_1$ and $W_2$, there exists a bijective correspondence
$K_{W_1W_2}$ which, given a $(W_1,W_2)$-compatible multiset $M$, produces a
pair $(P,Q)$ of bitableau of the same shape such that $P$ is a
$(\Omega_1,W_1)$-bitableau and $Q$ is a $(\Omega_2,W_2)$-bitableau.
A given value appears in the first (resp.~second) tableau exactly as many
times as it appears as the first (resp.~second) coordinate in $M$.
\end{thm}

\begin{proof}
This is essentially shown in \cite{BereleRemmel}; see also \cite{Remmel}.
While the references only deal with the case in which every element of
$\overline{W}$ is greater than every element of $W$, the proofs carry
over directly.
\end{proof}

\begin{rems}
We have switched the order of the tableaux usually used in the
Robinson-Schensted correspondence.
\end{rems}

\begin{rems}
The special cases $K_{\Omega_1\Omega_2}$ and $K_{\Omega_1\emptyset}$ are
known as the Knuth correspondence and the dual Knuth correspondence
\cite{Knuth:correspondence} respectively, while the special case
$K_{\emptyset\emptyset}$ is known as the Burge correspondence
\cite{Burge:correspondence,Gansner}.
\end{rems}

Since the above correspondence can be reduced to the Robinson-Schensted
correspondence (\cite{BereleRemmel}), all of the usual properties 
carry over.  We observe that the operations of inflation and
deflation carry over to bitableaux, and thus one can define
the dual $P^*$ of a bitableau.  That $P^*$ is indeed a bitableau
follows from the following:

\begin{thm}
Let $\iota_1:\Omega_1\to\Omega'_1$ and $\iota_2:\Omega_2\to\Omega'_2$ be
order-reversing maps; also, for a multiset $M$, define $M^t$ to be the
multiset in $\Omega_2\times \Omega_1$ obtained by exchanging the
coordinates.  For any finite multiset $M\subset\Omega_1\times\Omega_2$, the
following are equivalent:
\begin{align}
K_{W_1W_2}(M)&=(P,Q)\\
K_{W_1W_2}((\iota_1\times \iota_2)(M))&=(\iota_1(P^*),\iota_2(Q^*))\\
K_{W_2W_1}(M^t)&=(Q,P)\\
K_{W_2W_1}((\iota_2\times \iota_1)(M^t))&=(\iota_2(Q^*),\iota_1(P^*))\\
K_{\overline{W_1}\overline{W_2}}((\iota_1\times 1)(M))
&=(\iota_1(P^*)^t,Q^t))\\
K_{\overline{W_1}\overline{W_2}}((1\times \iota_2)(M))
&=(P^t,\iota_2((Q^*)^t))\\
K_{\overline{W_2}\overline{W_1}}((\iota_2\times 1)(M^t))
&=(\iota_2((Q^*)^t),P^t)\\
K_{\overline{W_2}\overline{W_1}}((1\times \iota_1)(M^t))
&=(Q^t,\iota_1((P^*)^t))
\end{align}
\end{thm}

\begin{rem}
For the cases $K_{\Omega_1\Omega_2}$ and $K_{\emptyset\emptyset}$, this
theorem appears in \cite{Gansner} and \cite{VoWhitney}.
\end{rem}

When $W_1=W_2=W$ and $M=M^t$, the two bitableau are the same; as in the
Robinson-Schensted correspondence, we can describe the number of odd-length columns:

\begin{thm}
Let $M$ be a finite multiset in $\Omega\times\Omega$ with $M=M^t$.  If
$\lambda$ is the common shape of $K_{WW}(M)$, then $f(\lambda')$ is
equal to the sum of (1) the number of elements of $M$ of the form
$(x,x)$ with $x\in W$ and (2) the number of $x\in \overline{W}$
that appear an odd number of times in $M$.
\end{thm}

\begin{proof}
This was known for $K_{\Omega_1\Omega_2}$ (\cite{Knuth:correspondence}) and
for $K_{\emptyset\emptyset}$ (\cite{Burge:correspondence}); it thus
follows in general.
\end{proof}

Let $M$ be a $(W_1,W_2)$-compatible multiset from $\Omega_1\times
\Omega_2$.  We observe that $l_{W_1W_2}(M)\le l$ whenever $M$ can be written as
the union of $l$ $(\overline{W_1},\overline{W_2})$-decreasing sequences.
(Again, we use $(W_1,W_2)$-compatibility).  This motivates the notation
$l^{(k)}_{W_1W_2}(M)$ for the size of the largest submultiset $M'$ of $M$
with $l^-_{\overline{W_1}\overline{W_2}}(M')\le k$, and similarly
for $l^{-(k)}_{\overline{W_1}\overline{W_2}}(M)$.

\begin{thm}
Let $M$ be a finite $(W_1,W_2)$-compatible multiset from
$\Omega_1\times\Omega_2$.  If $\lambda$ is the common shape of
$K_{W_1W_2}(M)$, then
\begin{align}
\sum_{1\le i\le k} \lambda_i = l^{(k)}_{W_1W_2}(M),\quad
\sum_{1\le i\le k} \lambda'_i = l^{-(k)}_{\overline{W_1}\overline{W_2}}(M).
\end{align}
\end{thm}

\begin{rem}
For the Robinson-Schensted correspondence, this was proved in \cite{Greene};
the extension to the general case is in \cite{BereleRemmel}.
\end{rem}

In particular, $\lambda_1$ gives the length of the longest
$(W_1,W_2)$-increasing sequence.

\begin{proof}[Proof of Theorem \ref{thm:kurtjg}]
We generalize the argument of \cite{kurtj:shape}.  Consider, for instance,
the case $\tsymmO$.  In this case, the probability that our random multiset
$M(q;\alpha)$ is equal to a given fixed multiset $M$ is
\[
\Prob(M(q;\alpha)=M)=
\Prob(M(q;\alpha)=\emptyset) \alpha^{f(\lambda')} q^P,
\]
where $K_{WW}(M)=(P,P)$ with $P$ of shape $\lambda$.  Thus
\begin{align}
\Prob(l^\symmO_W(q;\alpha)\le l)
&=
\sum_{\substack{\lambda_1\le l\\
\lambda_{WW}(M)=\lambda}}
\Prob(M(q;\alpha)=M)\\
&=
\Prob(M(q;\alpha)=\emptyset)
\sum_{\substack{\ell(\lambda)\le l\\
P\in B_{\lambda'}(\Z^+,W)}}
\alpha^{f(\lambda)} q^P\\
&=
Z^\symmO_W(q;\alpha)
\sum_{\ell(\lambda)\le l}
\alpha^{f(\lambda)}
s_{\lambda}(q_{\overline{W}}/q_W).
\end{align}

Similar arguments hold for types $\tsymmU$ and $\tsymmS$.  For $\tsymmUU$ and
$\tsymmu$, we also need the fact that $\tilde{s}_\lambda(x/y)$ can be
defined in terms of a sum over self-dual bitableaux.  Clearly, it
suffices to show that the usual correspondence between self-dual tableaux
and pairs of tableaux extends to bitableaux.

Define a domino $(\Omega,W)$-bitableau of shape $\lambda$ to be a tiling of
the diagram of $\lambda$ with labelled dominos such that the labels
increase weakly in each row and column, and such that (1) for $x\in W$, no
column hits more than one domino labelled $x$, and (2) for $x\notin W$, no
row hits more than one domino labelled $x$.  We readily verify (using the
fact that deflation preserves the bitableau property) that we have a
bijection between self-dual bitableaux and domino bitableaux.

To proceed from domino bitableaux to pairs of ordinary bitableaux, it
suffices to show that the usual correspondence preserves the bitableau
property.  But this follows from the fact that the correspondence is valid
for column-strict tableaux (and thus for row-strict tableaux by symmetry),
as remarked in \cite{StantonWhite}.

Combining these bijections, we obtain the desired correspondence,
and thus the theorem for $\tsymmUU$ and $\tsymmu$.
\end{proof}

Theorem \ref{thm:kurtjg} motivates the following change of notation:
\begin{align}
l^\symmU_{WW'}(q;q') &\to
l^\symmU(q_{\overline{W}}/q_{W};q'_{\overline{W'}}/q'_{W'})\\
l^\symmUU_{WW'}(q;q') &\to
l^\symmUU(q_{\overline{W}}/q_{W};q'_{\overline{W'}}/q'_{W'})\\
l^\symmO_W(q;\alpha) &\to
l^\symmO(q_{\overline{W}}/q_{W};\alpha)\\
l^\symmS_W(q;\beta) &\to
l^\symmS(q_{\overline{W}}/q_{W};\beta)\\
l^\symmu_W(q;\alpha,\beta) &\to
l^\symmu(q_{\overline{W}}/q_{W};\alpha,\beta)
\end{align}
It is somewhat startling that, even though the new notation in principle
discards information, the resulting distributions are in fact exactly the
same.  For the involution cases, the integral formulae of Section
\ref{sec:symmfunc} imply further that:

\begin{cor}\label{cor:kurtjg}
For any valid parameter choices, the following pairs of random variables
have the same distribution:
{
\allowdisplaybreaks
\begin{align}
l^\symmO(q/r;\alpha)
&\sim
l^\symmO(q/\alpha,r;0),\\
\lfloor \frac{1}{2} l^\symmS(q/r;\beta)\rfloor
&\sim
\frac{1}{2} l^\symmS(q/r;0),\\
\lceil \frac{1}{2} l^\symmS(q/r;\beta)\rceil
&\sim
\frac{1}{2} l^\symmS(\beta,q/r;0),\\
\lfloor \frac{1}{2} l^\symmu(q/r;\alpha,\beta)\rfloor
&\sim
l^\symmU(q/\alpha,r;q/r)\\
\lceil \frac{1}{2} l^\symmu(q/r;\alpha,\beta)\rceil
&\sim
l^\symmU(\beta,q/\alpha,r;q/r)
\end{align}
}
where, for instance $l^\symmO(\alpha,q/r;0)$ corresponds to a process in
which $\alpha$ has been inserted into the sequence $q$.
\end{cor}

\begin{proof}
We have:
\begin{align}
\Pr(l^\symmO(q/r;\alpha)\le l)
&=
Z^\symmO(q/r;\alpha)
\Exp_{U\in O(l)} \det((1+\alpha U)H(U;q/r))\\
&=
Z^\symmO(q/\alpha,r;0)
\Exp_{U\in O(l)} \det(H(U;q/\alpha,r))\\
\end{align}
and similarly for the other cases.
\end{proof}

\begin{rems}
By the arguments of Theorem \ref{thm:fixinv}, one can prove
something stronger.  Namely, the joint distribution of the lengths of the
odd-numbered rows is the same for $\tsymmO(q/r;\alpha)$ and for
$\tsymmO(q/\alpha,r;0)$, and similarly for the other cases.
For the even-numbered rows, the argument of Theorem \ref{thm:sfeig2inv}
shows that the joint distribution is independent of $\alpha$.
In the Laguerre limit (see below), we obtain the following fact.
Consider the ``matrix ensemble'' with joint eigenvalue density
(on $[0,\infty)$) proportional to
\[
\pf(\sgn(x_k-x_j)e^{A|x_k-x_j|})_{1\le j,k\le 2n}
\prod_{1\le j<k\le 2n} (x_k-x_j)
\prod_{1\le j\le 2n} e^{-C x_j},
\]
where $A$ and $C$ are parameters with $C>\max(A,0)$; note that if
$x_1<x_2<\ldots x_{2n}$, then
\[
\pf(\sgn(x_k-x_j)e^{A|x_k-x_j|})_{1\le j,k\le 2n}
=
\exp(A\sum_{1\le j\le 2n} (-1)^j x_{2j}).
\]
Then the joint distribution of the second, fourth, sixth, etc. largest
eigenvalues is independent of $A$.  Since for $A=0$ this ensemble is the
Laguerre orthogonal ensemble (LOE), while in the limit $A\to -\infty$, it
becomes the Laguerre symplectic ensemble (LSE), we find in particular that
the distribution of the second-largest eigenvalue of LOE is the same as the
distribution of the largest eigenvalue of LSE (since every eigenvalue of
LSE occurs twice).  For an alternate proof, and generalizations, see
\cite{ForRai99}.
\end{rems}
	
\begin{rems}
We can recover Theorems \ref{thm:trmoments} and \ref{thm:fixedpoissons}
from Theorem \ref{thm:kurtjg} by taking suitable limits.  For instance, for
$\tsymmU$, we take $q_i=q'_i=t/N$ for $1\le i\le N$, and $q_i=q'_i=0$
otherwise.  As $N\to\infty$, the resulting point process converges to the
usual Poisson process.  In Corollary \ref{cor:kurtjg}, if we take the
corresponding limit for the right-hand-sides, we obtain Poisson processes
in which, instead of adding extra diagonal points, we add extra \emph{side}
points.  Thus we obtain the fact that these two Poisson processes have
\emph{exactly} the same distribution.  For instance, if $n$ and $m$
are nonnegative integers, the distribution of the longest weakly increasing
subsequence is the same if we choose either: Pick $n$ points at random in
the triangle $0\le y\le x\le 1$, and $m$ points at random with $0\le x=y\le
1$; or: Pick $n$ points at random in the triangle $0\le y\le x\le 1$, and
$m$ points at random with $y=0$ and $0\le x\le 1$.  It is not at all clear
why these distributions should be the same.
\end{rems}

\begin{rems}
By taking a Poisson limit for only some of the variables, we obtain
increasing subsequence interpretations for the image of the Schur function
sums under homomorphisms of the form
\[
H(t;z) \mapsto e^{at} H(t;x) E(t;y)
\]
with $a$, $x_i$ and $y_i$ nonnegative (and satisfying the appropriate
additional conditions).  As this is the most general case in which the
images of $s_\lambda$ are all nonnegative (\cite{MacdonaldT3}),
this is presumably the most general case for which such an interpretation
can be given.
\end{rems}

\begin{rems}
In addition to the Poisson limit, another natural limit is the Laguerre
limit (\cite{kurtj:shape}), in which the random multiplicities are
exponentially distributed.  More precisely, one fixes integers $N$ (and
$N'$ as necessary), and considers a process in which $q$ and $q'$ are
constant on the respective intervals $[1,N]$ and $[1,N']$, and 0
otherwise, then takes the limit $q$, $q'\to 1$.  In the case
$\alpha=\beta=0$, one finds the following curious fact, analogous to
Theorem \ref{thm:trmoments}.

Given a complex $N\times N'$ matrix $M$, we define two decreasing sequences
of nonnegative real numbers.  $\Sigma(M)_i$ is the $i$th largest eigenvalue
of $M M^\dagger$, or 0 if $i>\min(N,N')$, while $\Delta(M)_i$ is
defined so that
\[
\sum_{1\le i\le j} \Delta(M)_i = \max_S \sum_{(k,l)\in S} |M_{kl}|^2,
\]
where $S$ ranges over unions of $j$ decreasing paths.

\begin{thm}
Let $N$ be a positive even integer.
Map $H$ (recall Section \ref{sec:symmiseq}) into transformations of
$N\times N$ complex matrices as follows:
\[
\diagup \mapsto (M\mapsto -M^t)\quad\text{and}\quad
\diagdown \mapsto (M\mapsto J M^t J)
\]
Let $G^\symmg(N)$ be the Gaussian distribution on matrices with the
appropriate symmetry.  Then for each $\symmg$, the distributions
$\Sigma(G^\symmg(N))$ and $\Delta(G^\symmg(N))$ are the same.
\end{thm}

\begin{proof} (Sketch)
That we can compute the distribution of $\Delta(G^\symmg(N))$ (as well as
the significance of this result) follows from the fact that it is the
Laguerre limit of the appropriate discrete process with symmetry $\tsymmg$.
Thus, as in \cite{kurtj:shape} for $\tsymmU$, we find that the appropriate
sum of pfaffians (see Section \ref{sec:pfaffians}) tends to a Riemann
integral.

On the other hand, the distributions $G^\tsymmg(N)$ for $\tsymmU$, $\tsymmO$,
and $\tsymmS$ are well-studied (these are unconstrained, antisymmetric, and
symmetric complex matrices, respectively), and in particular the
distribution $\Sigma(G^\symmg(N))$ is known in these cases (see
\cite{ForRai99}).  For $\tsymmUU$ and $\tsymmu$, we can perform simple row
and column operations (not changing $\Sigma(M)$) to reduce to the case of
$\tsymmU$.

We find that the two density functions we obtain are the same, proving the
theorem.
\end{proof}

\begin{rem}
For $\tsymmU$ and $\tsymmO$, $N$ need not be even, while
for $\tsymmU$ and $\tsymmUU$, the matrices need not be square.
\end{rem}
\end{rems}

\section{Invariants and increasing subsequences}\label{sec:invariants}

Consider the integral
\[
\Exp_{U\in U(l)} |\Tr(U)|^{2n}.
\label{eq:invdim}
\]
The function $|\Tr(U)|^{2n}$ is the character of a representation of
$U(l)$, since
\[
|\Tr(U)|^{2n}
=
\Tr(U^{\otimes n} \otimes \overline{U}^{\otimes n}),
\]
where we write $A^{\otimes n}$ for the tensor product of $A$ with itself
$n$ times.  It follows that \eqref{eq:invdim} (and thus $f^\symmU_{nl}$)
gives the dimension of the fixed subspace of that representation.
Equivalently, this is the dimension of the space $C_n(U(l))$ of operators on
$(\C^l)^{\otimes n}$ that commute with $U^{\otimes n}$ for all $U\in
U(l)$.

Given a permutation $\pi\in S_n$, we can associate an operator
$T_l(\pi)$ on $(\C^l)^{\otimes n}$ as follows:
\[
T_l(\pi) (v_1\otimes v_2\otimes \dots v_n)
=
(v_{\pi(1)}\otimes v_{\pi(2)}\otimes \dots v_{\pi(n)}).
\]
This operator clearly commutes with $U^{\otimes n}$ for $U\in U(l)$,
and thus $T_l$ extends to a map from $\C[S_n]$ to $C_n(U(l))$.
Indeed, this map is known to be surjective (see, e.g., \cite{Brauer});
i.e., the operators $T_l(\pi)$ span $C_n(U(l))$.  In general, however,
it is not injective.  For any subset $S\subset \{1,2,\ldots
n\}$, define two elements of $\C[S_n]$:
\begin{align}
E_S&\defeq \sum_{\substack{\pi\in S_n\\
\pi(x)=x,\ \forall x\notin S}} \sigma(\pi) \pi\\
H_S&\defeq \sum_{\substack{\pi\in S_n\\
\pi(x)=x,\ \forall x\notin S}} \pi.
\end{align}

\begin{lem}\label{lem:invkerU}
For $l\ge n$, $T_l$ is injective on $\C[S_n]$, while for $l<n$, the kernel
of $T_l$ contains all elements of the form $\pi E_S$ with $|S|>l$.
\end{lem}

\begin{proof}
That $T_l$ is injective for $l\ge n$ is straightforward; we simply observe
that if $v_1$, $v_2$,\dots $v_n$ are linearly independent vectors, then
the vectors
\[
T_l(\pi) (v_1\otimes v_2\otimes \dots v_n)
=
v_{\pi(1)}\otimes v_{\pi(2)}\otimes\dots v_{\pi(n)}
\]
are linearly independent, as $\pi$ ranges over $S_n$.

Thus, suppose $l<n$.  That $\pi E_S\in\ker T_l$ follows from the fact
that any tensor product of $|S|>l$ basis vectors must contain at least
one basis vector more than once, so will be taken to 0 by $E_S$.
\end{proof}

\begin{rem}
It follows from the proof of Theorem \ref{thm:invU} below that these
elements also span the kernel.
\end{rem}

Using this fact, we obtain:

\begin{thm}\label{thm:invU}
For any nonnegative integers $l$ and $n$, the set
\[
\{T_l(\pi):\pi\in S_n|\ell^-(\pi)\le l\}
\]
is a basis of $C_n(U(l))$.
\end{thm}

\begin{proof}
We first need to show that, given any permutation $\pi$ with a long
decreasing subsequence, we can express $T_l(\pi)$ as a linear combination
of $T_l(\pi')$ with $\pi'$ ranging over permutations without long
decreasing subsequences.

Let $\pi$ be such a permutation, and let $S$ be a subset of size $l+1$
on which $\pi$ is decreasing.  By the lemma, it follows that
\begin{align}
T_l(\pi) &= T_l(\pi-\pi E_S)\\
&=
-\sigma(\pi)
\sum_{\substack{\pi'\in S_n\\
\pi'(x)=\pi(x),\ \forall x\notin S\\\pi'\ne \pi}} \sigma(\pi') T_l(\pi').
\end{align}
Now each permutation $\pi'$ on the right hand side agrees with $\pi$
outside $S$, but is no longer decreasing on $S$.  It follows that each
$\pi'$ has strictly fewer inversions than $\pi$ (reduce to the case in which
$\pi'$ differs from $\pi$ by a 2-cycle).

It follows that if we iterate this reduction, we will eventually
obtain a linear combination of permutations that cannot be reduced.
But this is precisely the desired result.

It remains only to show that the elements $T_l(\pi)$ are linearly
independent, which we will do via a triangularity argument.  If we choose a
basis of $V$, we can view $T_l$ as the restriction to $\{1,2,\ldots l\}^n$
of the corresponding action of $S_n$ on the set $W_n=\N^n$.  Now, to a
given permutation $\pi$, we associate two words $w_1(\pi)$ and $w_2(\pi)$
as follows: $w_1(\pi)_j$ is equal to the length of the longest decreasing
subsequence of $\pi$ starting with $j$; similarly $w_2(\pi)_j$ is equal to
the length of the longest decreasing subsequence of $\pi$ starting in
position $j$.  (Since $\pi$ has longest decreasing subsequence of length $\le l$, these are indeed in $\{1,2,\ldots l\}^n$.)
We easily see that $w_2(\pi) = \pi(w_1(\pi))$, and
thus that $T_l(\pi)$ has coefficient $1$ on the pair $(w_1(\pi),w_2(\pi))$.
The claim is then that any other permutation taking $w_1(\pi)$ to
$w_2(\pi)$ has strictly more inversions than $\pi$.

Define the number of inversions $i(w)$ of a word $w$ to be the number of
coordinate positions $i<j$ such that $w_j<w_i$.  By standard arguments, we
find that if $\pi(v)=w$, then $i(w)\le i(v)+i(\pi)$; equality holds if for
any pair $i<j$ such that $\pi_i>\pi_j$, we have $v_i<v_j$ and $w_i>w_j$.
We immediately deduce that (a) $i(w_2(\pi))=i(w_1(\pi))+i(\pi)$ and
(b) if $\pi'(w_1(\pi))=w_2(\pi)$, then $i(\pi')\ge i(\pi)$.
Finally, if $\pi'(w_1(\pi))=w_2(\pi)$ with $i(\pi')=i(\pi)$,
then $\pi'$ has not only the same number of inversions as $\pi$,
but indeed the same {\emph set} of inversions; it follows that
$\pi'=\pi$.
\end{proof}

Probably the most important consequence of Theorem \ref{thm:invU} is
that it gives an elementary proof that
\begin{cor}\label{cor:invUc}
For any integers $n$ and $l$,
\[
E_{U\in U(l)} |\Tr(U)|^{2n} = \dim(C_n(U(l))) = f_{nl}.
\]
\end{cor}

\begin{rems}
A reduction algorithm closely related to that used above appeared in
\cite{Re78}, for an application to P.I. algebras.  The connection
to invariant theory, as well as the various generalizations given below
appear to be new, however.
\end{rems}

\begin{rems}
A different basis for $l=2$ (based on one of the many combinatorial
interpretations of the Catalan numbers) was given in \cite{invarsl2}.
\end{rems}

\begin{rems}
We could, of course, just as easily have used the permutations without long
{\it increasing} subsequences to form the basis.  The current choice
has the merit of giving a basis containing the identity and closed under
taking inverses, as well as making the proof of linear independence
somewhat cleaner.
\end{rems}

\begin{rems}
While we defined everything over $\C$, we observe that both
the representation $T_l$ and the above reduction algorithm are
actually defined over $\Z$.
\end{rems}

\begin{rems}
There are also, of course, analogues of Corollary \ref{cor:invUc}
associated with Theorems \ref{thm:invUm} through \ref{thm:invSm};
we leave the details to the reader.
\end{rems}

More generally, Theorem \ref{thm:kurtjg} implies
\[
Z^\symmU_{WW'}(q;q')^{-1} \Prob(l^\symmU_{WW'}(q;q')\le l)
=
E_{U\in U(l)}
\det(H(U;q_{\overline{W}}/q_W)
H(U^\dagger;q'_{\overline{W'}}/q'_{W'})).
\]
Consider the coefficient of a monomial
\[
\prod_i q_i^{\nu_i} (q'_i)^{\nu'_i}
\]
in the right-hand-side.  By the properties of $H(t;x/y)$, this is
\[
E_{U\in U(l)}
\prod_{i\in W} e_{\nu_i}(U)
\prod_{i\in \overline{W}} h_{\nu_i}(U)
\prod_{i\in W'} \overline{e_{\nu'_i}(U)}
\prod_{i\in \overline{W'}} \overline{h_{\nu'_i}(U)}.
\]
Again, this is the expectation of a character, and thus computes the
dimension of a space of invariants.  The corresponding coefficient of the
left-hand-side counts $(W,W')$-compatible multisets without
$(W,W')$-increasing subsequences of length $l+1$, in which $i$ appears
$\nu_i$ times as a first coordinate and $\nu'_i$ times as a second
coordinate.  Thus to obtain the analogue of Theorem \ref{thm:invU}, we
first need an analogue of $T_l$ for multisets.

Given a composition $\nu$, we define a partition $S(\nu)$ of
$\{1,2,\ldots |\nu|\}$ by
\[
S(\nu)_i = \{j:\sum_{k<i} \nu_k<j\le \sum_{k\le i} \nu_k\}.
\]
Associated to this partition is an operator
\[
\Pi(\nu) = \prod_{i \in W} E_{S(\nu)_i} \prod_{i\in \overline{W}}
H_{S(\nu)_i}
\]
Then the representation of $U(l)$ in question is simply the action by
conjugation of $U^{\otimes n}$ on operators of the form $T_l(\Pi(\nu')) A
T_l(\Pi(\nu))$.  In particular, the invariant subspace is spanned by
operators of the form $T_l(\Pi(\nu')\pi\Pi(\nu))$ with $\pi\in S_n$. We
note that if two elements of $S_i$ with $i\in W$ each map to elements of
the same $S'_j$ with $j\not\in W'$, then $\Pi(\nu') \pi \Pi(\nu)=0$; and
similarly if $i\not\in W$ and $j\in W'$.  It follows that operators of the
form $\Pi(\nu')\pi\Pi(\nu)$ are, up to sign, in one-to-one correspondence
with $(W,W')$-compatible multisets with content $(\nu,\nu')$.  In
particular, given a $(W,W')$-compatible multiset $M$ with content
$(\nu,\nu')$, we obtain a element of $\C[S_n]$ which we denote $T_l(M)$.
Writing $\calM^\symmU_{WW'}(\nu;\nu')$ for the set of such multisets,

\begin{thm}\label{thm:invUm}
For any nonnegative integers $l$ and compositions $\nu$ and $\nu'$,
\[
\{T_l(M):M\in \calM^\symmU_{WW'}(\nu;\nu')|\ell^-_{WW'}(M)\le l\}
\]
is a basis of $\Pi(\nu') C_n(U(l)) \Pi(\nu)$.
\end{thm}

\begin{proof}
Given a multiset $M$ in $\Z^+\times \Z^+$, we define an inversion of $M$
to be a pair of elements $(x,y)\in M$, $(z,w)\in M$ with $x<z$ and $y>w$
(i.e., a strictly decreasing subsequence of length 2).

Now, suppose $M\in \calM_{WW'}(n;n')$ has a $(W,W')$-decreasing subsequence
of length $l+1$.  Choose a permutation $\pi$ corresponding to $M$,
and let $S\subset\{1,2,\ldots n\}$ be the set of positions of $\pi$
corresponding to the $(W,W')$-decreasing subsequence.  As before,
we have
\[
T(\Pi(\nu')\pi E_S\Pi(\nu)) = 0.
\]
Now, we readily see that for any permutation that appears in the
left-hand-side, the corresponding multiset has at most as many inversions
as $M$.  Indeed, the number of inversions is strictly smaller
unless the corresponding multiset is \emph{equal} to $M$.  Thus it
remains only to show that the terms corresponding to $M$ do not
cancel.  The only way that the term
\[
\sigma(\pi') \Pi(\nu') \pi \pi' \Pi(\nu)
\]
can correspond to $M$ is if we can write
\[
\pi' = (\prod_i \pi^{-1} \pi'_{1i} \pi) (\prod_i \pi'_{2i}),
\]
where $\pi'_{2i}$ fixes the complement of $S\cap S(\nu)_i$ and
$\pi'_{1i}$ fixes the complement of $\pi(S)\cap S(\nu')_i$.  Now, by
the definition of $(W,W')$-increasing subsequence, $S\cap S(\nu)_i$
contains at most one element when $i\notin W$; similarly
$\pi(S)\cap S(\nu')_i$ contains at most one element when $i\notin W'$.
In other words, we can write
\[
\pi' = (\prod_{i\in W'} \pi^{-1} \pi'_{1i} \pi) (\prod_{i\in W} \pi'_{2i}),
\]
But then
\[
\sigma(\pi') \Pi' \pi \pi' \Pi = \Pi' \pi \Pi,
\]
as desired.

The proof of linear independence is analogous to that in Theorem
\ref{thm:invU}.  Of the permutations associated to $M$, there is a
unique one ($\pi(M)$) such that each element of $W$ and $W'$ induces a
decreasing subsequence and each element of $\overline{W}$ and
$\overline{W'}$ induces an increasing subsequence.  We define
$w_1(M)=w_1(\pi(M))$, $w_2(M)=w_2(\pi(M))$, and observe that any other
multiset $M'$ with a nonzero coefficient at $(w_1(M),w_2(M))$
satisfies $i(\pi(M'))>i(\pi(M))$.
\end{proof}

%
%
%
\begin{rems}
It is possible to renormalize the basic invariants $T_l(M)$ in such a way
that the reduction algorithm is integral.  We need simply divide $T_l(M)$
by
\[
\prod_{i\in W,j\in W'} |S(\nu)_i\cap \pi(S(\nu')_j)|!,
\]
where $T_l(M)=\Pi(\nu')\pi\Pi(\nu)$.
\end{rems}

\begin{rems}
The space $\Pi(\nu') C_n(U(l)) \Pi(\nu)$ can be thought of as the space of
simultaneous (multilinear) invariants of a collection of symmetric and
antisymmetric tensors, some covariant and some contravariant.  Of special
interest is the case in which $\nu'$ is the composition $l^k$, with $1$,
$2$,\dots $k\in W'$.  In this case, the space
\[
(E_{\{1,2,\ldots l\}} E_{\{l+1,l+2,\ldots 2l\}}\dots) C_n(U(l)) \Pi(\nu)
\]
corresponds to \emph{relative} invariants of a collection of symmetric and
antisymmetric covariant tensors; i.e., transforming the tensors multiplies
the invariant by a power of the determinant.  There is a known
algorithm (the straightening algorithm \cite{RotaSturmfels}) for
computing a basis of such invariants.  In fact, we observe that the
resulting basis is, up to constant factors, the same as our basis.
Thus our algorithm can be viewed as a generalization of this algorithm
(different from the generalization to the ``fourfold'' algebra of
\cite{Rotabook}). Similarly, the algorithms below for the orthogonal and
symplectic groups can be thought of as straightening algorithms for those
groups.
\end{rems}

\begin{rems}
There is also a ``quantum'' analogue of this result.  If one replaces the
unitary group $U(l)$ by the quantum enveloping algebra $U_q({\mathfrak gl}_l)$,
the role of the symmetric group is now played by the Hecke algebra
\cite{Jimbo}.  As long as $S$ consists of consecutive elements, there is no
difficulty in defining $E_S$ and $H_S$ (these are idempotents corresponding
to 1-dimensional characters of parabolic subalgebras).  We find that the
kernel of the quantum $T_l$ is the ideal generated by $E_S$ with
$S=\{k,k+2,\ldots k+l\}$, $1\le k\le n-l$.  Using the appropriate
normalization, we obtain a reduction algorithm integral over $\Z[q]$.
(We also have linear independence whenever $q$ is not a root of unity.)
Of particular interest is the ``crystal limit'' $q=0$.  Under that
specialization, the relations take the form
$T_l(M) = 0$,
with $M$ ranging over multisets with long decreasing subsequences.  This
case surely merits further investigation, given the connection
\cite{quantumRSK} between the crystal limit of the quantum straightening
algorithm and the Robinson-Schensted-Knuth correspondence.  It would
also be interesting to understand the analogues for the quantum orthogonal
and symplectic groups.
\end{rems}

For the orthogonal group $O(l)$, there is a ``basic'' invariant
$T^\symmO_l(\tau)$ associated to any fixed-point-free involution $\tau$ in
$S_{2n}$ such that the basic invariants span the invariant space of
$U^{\otimes 2n}$ (which we denote by $F_{2n}(O(l))$.  These transform under
$T_l(\pi)$ as:
\[
T_l(\pi) T^\symmO_l(\tau) = T^\symmO_l(\pi^{-1} \tau\pi).
\]

We write $\pi\cdot \tau$ for the corresponding action on the formal span of
fixed-point-free involutions.  For any composition $\nu$ and nonnegative
integer $a$, we define
\[
\calM^\symmO_W(\nu;a)
\]
to be the set of symmetric $(W,W)$-compatible multisets with content
composition $\nu$ and with $a$ ``fixed points'' (i.e., $(i,i)\in M$ such
that $i\in W$, together with $(i,i)$ with $i\notin W$ such that $(i,i)$ has
odd multiplicity.)  Clearly, there is a correspondence (up to sign) between
$\calM^\symmO_W(n_i;0)$ and elements of the form $\Pi\cdot \tau$.

\begin{thm}\label{thm:invO}
For any nonnegative integer $l$ and any composition $\nu$,
\[
\{T^\symmO_l(M):
M\in \calM^\symmO_{W}(\nu;0)|\ell_{WW}(M)\le l\}
\]
is a basis of $\Pi(\nu) F_{2n}(O(l))$.
\end{thm}

\begin{proof}
The argument is analogous.  In eliminating a given increasing subsequence,
we replace both it and its reflection through the diagonal by
non-increasing subsequences, so the number of inversions increases.

The proof of linear independence is analogous to that of Theorem
\ref{thm:invUm}.  We simply switch the reverse the inequalities on the
second coordinates before applying the above arguments.
\end{proof}

This corresponds to taking the coefficient of a monomial $q^\nu$ in
the identity
\[
Z^\symmO_W(q;\alpha)^{-1}
\Prob(l^\symmO_W(q;\alpha)\le l)
=
E_{U\in O(l)} \det(H(U;q_{\overline{W}}/\alpha,q_W))
\]
To handle the general case, we need to map an element of
$M^\symmO_W(\nu;a)$ (corresponding to the monomial $q^\nu \alpha^a$) to an
element of
\[
\Pi^\symmO(\nu;a) F_{2n}(O(l)),
\]
where we define
\[
\Pi^\symmO(\nu;a) = \Pi(\nu) E_{\{|\nu|+1,|\nu|+2,\ldots |\nu|+a\}}.
\]
But this is simple: simply convert each fixed point $(i,i)$ to a pair of
elements $(i,\infty)$ and $(\infty,i)$, and apply $T^\symmO_l$.

\begin{thm}\label{thm:invOm}
For any nonnegative integers $l$ and $a$ and any composition $\nu$,
\[
\{T^\symmO_l(M):
M\in \calM^\symmO_W(\nu;a)|\ell_{WW}(M)\le l\}
\]
is a basis of $\Pi^\symmO(\nu;a) F_{2n}(O(l))$.
\end{thm}

\begin{proof}
The main difficulty here is that the na\"\i{}ve extension of the above
algorithm is no longer guaranteed to terminate; for instance, 
corresponding to the increasing subsequence $13$ of $132$, we have
the identity
\[
T^\symmO_2(132)=T^\symmO_2(213).
\]
Since $13$ is an increasing subsequence of both sides, we could clearly
loop indefinitely.

The solution is to require that the increasing subsequence consist entirely
of points $(x,y)$ with $x\le y$.  Even with this restriction, the above
proof does not entirely carry over; for instance, both $132$ and $213$ have
the same number of inversions (i.e., 1).

For an involution $\tau$, denote by $S_\le(\tau)$ the set of $i$ with $i\le
\tau(i)$.  Given a multiset $M$, we then define $M_\le(M)$ to be the
multiset corresponding to $S_\le(\tau)$ where $\tau$ corresponds to $M$.
Given two multisets $M_1$ and $M_2$ of the same size on the same totally
ordered set, we write $M_1\le M_2$ to indicate that we can identify
elements of $M_1$ and $M_2$ in such a way that each element of $M_1$ is
$\le$ the corresponding element of $M_2$.  Then the theorem follows from
the following observation:

If $M'$ is one of the multisets obtained after eliminating an increasing
subsequence of $M$, then $M_\le(M')\le M_\le(M)$.  If equality occurs, then
either $M'=M$ or $M'$ has strictly more inversions than $M$.

Linear independence follows as above.
\end{proof}

For the symplectic group $Sp(2l)$, the basic invariants again correspond to
involutions, but now transform under $T_l(\pi)$ as:
\[
T_{2l}(\pi) T^\symmS_l(\tau) = \pm T^\symmS_{2l}(\pi^{-1} \tau\pi),
\]
with
\[
T_{2l}(\pi) T^\symmS_l(\tau) = \sigma(\pi) T^\symmS_l(\tau)
\]
whenever $\pi$ commutes with $\tau$.  For any composition $\nu$ and
nonnegative integer $b$, we define
\[
\calM^\symmS_W(\nu;b)
\]
to be the set of symmetric $(\overline{W},\overline{W})$-compatible
multisets with content composition $\nu$ and with $b$ fixed points.

\begin{thm}\label{thm:invS}
For any nonnegative integer $l$ and composition $\nu$,
\[
\{T^\symmS_l(M):
M\in \calM^\symmS_{W}(\nu;0)|\ell^-_{WW}(M)\le 2l\}
\]
is a basis of $\Pi(\nu) F_{2n}(Sp(2l))$.
\end{thm}

\begin{proof}
Analogous.  The only issue is that the long decreasing subsequence we
eliminate must be symmetric about the diagonal (clearly always possible).

For linear independence, we choose a symplectic basis of $V$ indexed
$v_{\pm 1}$, $v_{\pm 2}$.  We thus find that the nonzero coefficients of an
involution $\tau$ correspond to words of length $2n$ on $\pm \Z^+$ such
that $\tau(w)=-w$.  The word $w(\tau)$ associated to an involution
is now defined so that $|\tau_j|$ is equal to half the length of the
longest symmetric decreasing subsequence starting or ending with $j$;
the sign is positive if the sequence ends with $j$, and negative
otherwise.  Other than this, the arguments are analogous.
\end{proof}

This extends easily to the case when fixed points are allowed; in this
case, we define
\[
\Pi^\symmS_{2l}(\nu;b)
=
\Pi(\nu) H_{|\nu|+1,|\nu|+2,\dots |\nu|+b}.
\]
When eliminating (symmetric) decreasing subsequences, if the subsequence
has even length, we can proceed as above; otherwise, we permute only those
elements not corresponding to $\infty$.  In either case, we convert an
decreasing subsequence to a non-decreasing subsequence.  (And the linear
independence argument carries over.)  Thus

\begin{thm}\label{thm:invSm}
For any nonnegative integers $l$ and $b$ and any composition $\nu$,
\[
\{T^\symmS_l(M):
M\in \calM^\symmS_W(\nu;b)|\ell^-_{WW}(M)\le 2l\}
\]
is a basis of $\Pi^\symmS(\nu;b) F_{2n}(Sp(2l))$.
\end{thm}

It is not entirely clear how to proceed for either of $\tsymmUU$ or
$\tsymmu$.  For $\tsymmUU$, or more generally $U(p)\times U(q)$, the
centralizer algebra corresponds as expected to a certain representation
$T_{p,q}$ of the hyperoctahedral group, in which permutations are mapped
via $T_{p+q}$, and sign changes correspond to an operator with eigenvalues
1 ($p$ times) and $-1$ ($q$ times).  Given a set $S\subset \{1,2,\ldots
n\}$, we can define two elements $E^{\pm}_S\in \C[H_n]$.  Each is a sum
over elements of $H_n$ that fix the complement of $S$; in $E^+$, we
multiply by the sign of the corresponding permutation, while in $E^-$, we
further multiply by the number of sign changes.  Since the group
$U(p)\times U(q)$ is a direct product, we obtain:

\begin{lem}
The kernel of $T_{p,q}$ on $\C[S_n]$ is spanned by the elements
$\pi E^+_S$ with $|S|>p$ and by the elements $\pi E^-_S$ with $|S|>q$.
\end{lem}

\begin{rem}
The point is that $\pi E^+_S$ can be used to reduce invariants of $U(p)$,
while $\pi E^-_S$ can be used to reduce invariants of $U(q)$.  Since
every invariant of $U(p)\times U(q)$ can be expressed in terms of the
invariants of the respective factors, we are done.
\end{rem}

Even in the cases of interest, however, ($q=p$ or $q=p+1$), it is not clear
how to use these relations to eliminate hyperoctahedral permutations with
long decreasing subsequences; in particular, the invariants do not span the
kernel of $T_{p,q}$ on $\Z[S_n]$ in general, but only on $\Z[1/2][S_n]$.
Similar remarks apply to $\tsymmu$.  As in the other involution cases, the
invariant space for $\tsymmu$ is associated to a (twisted) Gelfand pair in
$H_{2n}$; the relevant subgroup of $H_{2n}$ is the centralizer of an
element of $S^\symmu_n$.  (For $\tsymmO$ and $\tsymmS$, the Gelfand pair
is $(S_{2n},H_n)$, twisted in the $\tsymmS$ case.)

We observe that in each case, $\ker(T^\symmg_l)\subset
\ker(T^\symmg_{l-1})$ for all $l$, and $\ker(T^\symmg_n)=0$.  We thus
obtain the following theorem, which is a formal analogue of the Szeg\"o
limit theorem:

\begin{thm}\label{thm:formalSzego}
We have the following limits of formal power series:
\begin{align}
\lim_{l\to\infty}
E_{U\in U(l)} \det(H(U;x/y)H(U^\dagger;z/w))
&=
\prod_{j,k} (1-x_jz_k)^{-1}(1-y_jw_k)^{-1}
\prod_{j,k} (1+x_jw_k)(1+y_jz_k).\\
\lim_{l\to\infty}
E_{U\in O(l)} \det(H(U;x/y))
&=
\prod_{j,k} (1+x_jy_k)
\prod_{j\le k} (1-x_jx_k)^{-1}
\prod_{j<k} (1-y_jy_k)^{-1}.\\
\lim_{l\to\infty}
E_{U\in Sp(2l)} \det(H(U;x/y))
&=
\prod_{j,k} (1+x_jy_k)
\prod_{j<k} (1-x_jx_k)^{-1}
\prod_{j\le k} (1-y_jy_k)^{-1}.
\end{align}
More precisely, all coefficients of degree $\le 2l$ agree, and each limit
is monotonic in all coefficients.
\end{thm}

\begin{rems}
Under the homomorphism $p_j(x/y)\mapsto f_i$, $p_j(z/w)\mapsto g_i$,
we obtain
\begin{align}
\lim_{l\to\infty}
E_{U\in U(l)} 
\exp(\sum_j f_j \Tr(U^j)/j +\sum_j g_j \Tr(U^{-j})/j)
&=
\exp(\sum_j f_j g_j/j)\\
\lim_{l\to\infty}
E_{U\in O(l)} 
\exp(\sum_j f_j \Tr(U^j)/j)
&=
\exp(\sum_j (f^2_j+f_{2j})/2j)\\
\lim_{l\to\infty}
E_{U\in Sp(2l)}
\exp(\sum_j f_j \Tr(U^j)/j)
&=
\exp(\sum_j (f^2_j-f_{2j})/2j)
\end{align}
and again coefficients with (weighted) degree $\le 2l$ agree.  This fact
was proved via representation theory in \cite{DiaconisShahshahani} (note
that the degree bounds given there are incorrect); the connection to
Szeg\"o's limit theorem was observed in \cite{Jo1}.  The monotonicity
result is new.
\end{rems}

\begin{rems}
There are, of course, also hyperoctahedral analogues, both of which simply
reduce to the $\tsymmU$ case.
\end{rems}

We can also give a partial extension of the above results to the ``super''
analogues of the classical groups.  For the unitary supergroup
$U(l/k)$, the centralizer algebra $C_n(U(l/k))$ is again spanned by
permutations, under a particular representation $T_{l/k}$ of $S_n$.

\begin{thm}
The representations that appear in $T_{l/k}$ with positive multiplicity
are those indexed by partitions $\lambda$ with $\lambda_{l+1}<k+1$.
\end{thm}

We then have:

\begin{thm}
For any nonnegative integers $l$ and $k$ and compositions $\nu$ and $\nu'$,
the dimension of the space $\Pi(\nu') C_n(U(l/k))\Pi(\nu)$ is equal to the
number of multisets $M\in \calM^\symmU_{WW'}(\nu;\nu')$ such that
\[
\lambda_{WW'}(M)_{l+1}<k+1.
\]
\end{thm}

\begin{proof}
If $T_\lambda$ is the representation of $S_n$ corresponding to the
partition $\lambda$, then the dimension of
\[
T_\lambda(\Pi(\nu') S_n \Pi(\nu))
\]
is given by the number of pairs of bitableau of shape $\lambda$
with respective content $\nu$ and $\nu'$.  In other words, by the
generalized Knuth correspondence, this is equal to the number of multisets
$M\in \calM^\symmU_{WW'}(\nu;\nu')$ with
\[
\lambda_{WW'}(M)=\lambda.
\]
The result follows by summing over $\lambda$.
\end{proof}

\begin{rems}
The obvious conjecture that those multisets give a basis under
$T_{l/k}$ is not true in general (not even for $C_4(U(1/1))$);
thus Theorem \ref{thm:invUm} does not immediately extend to the
supergroup case.  It is also not clear whether there is a simple
description of the kernel of $T_{l/k}$.
\end{rems}

\begin{rems}
This theorem can be thought of as a formal statement along the lines of:
\[
E_{U\in U(l/k)} \det(H(U;x)H(U^\dagger;y))
=
\sum_{\lambda_{l+1}<k+1} s_\lambda(x)s_\lambda(y).
\]
It would be nice to make this statement precise.  It is also interesting to
speculate on the possibility of a ``super'' analogue of orthogonal
polynomials and Toeplitz determinants.
\end{rems}

Similarly, for the orthosymplectic supergroup $OSp(l/2k)$ (the full group,
not just the component of the identity), the invariant space
$F_n(OSp(l/2k))$ is the image of fixed-point-free involutions under a map
$T^\symmO_{l/2k}$ for which
\[
T_{l/2k}(\pi)T^\symmO_{l/2k}(\tau)
=
T^\symmO_{l/2k}(\pi^{-1} \tau \pi).
\]
This thus induces an action of $S_{2n}$ on $F_{2n}(O(l/2k))$, for which:

\begin{thm}
$F_{2n}(OSp(l/2k))$ splits into $S_n$-irreducible submodules as the direct sum
\[
\bigoplus_{\substack{\lambda\vdash n\\ \lambda_{l+1}\le k}} T_{2\lambda}.
\]
\end{thm}

It then follows that

\begin{thm}
For any nonnegative integers $l$ and $k$ and compositions $\nu$ and $\nu'$,
the dimension of the space $\Pi(\nu) F_{2n}(OSp(l/2k))$ is equal to the
number of multisets $M\in \calM^\symmO_{WW'}(\nu;0)$ such that
\[
\lambda_{WW'}(M)_{l+1}<2k+1.
\]
\end{thm}

In ``integral'' form, this reads:
\[
E_{U\in OSp(l/2k)} \det(H(U;x))
=
\sum_{\lambda_{l+1}\le k} s_{2\lambda}(x).
\]



\end{document}